\documentclass[12pt]{amsart}
\usepackage{amssymb,mathrsfs}
\title{The Weil representations of the Jacobi group }
\author{ Jae-Hyun Yang}
\address{Department of Mathematics, Inha University, Incheon 402-751, Korea}
\email{jhyang@inha.ac.kr }

\begin{document}

\thanks{This work was supported by Inha University Research Grant}


\maketitle
\begin{abstract} The Jacobi group is the semi-direct product of the symplectic group and
the Heisenberg group. The Jacobi group is an important object in the framework of
quantum mechanics, geometric quantization and optics. In this paper, we study the
Weil representations of the Jacobi group and their
properties. We also provide their applications to the theory of
automorphic forms on the Jacobi group and representation theory of
the Jacobi group.
\end{abstract}

\newtheorem{theorem}{Theorem}[section]
\newtheorem{lemma}{Lemma}[section]
\newtheorem{proposition}{Proposition}[section]
\newtheorem{remark}{Remark}[section]
\newtheorem{definition}{Definition}[section]

\renewcommand{\theequation}{\thesection.\arabic{equation}}
\renewcommand{\thetheorem}{\thesection.\arabic{theorem}}
\renewcommand{\thelemma}{\thesection.\arabic{lemma}}
\newcommand{\bbr}{\mathbb R}
\newcommand{\bbs}{\mathbb S}
\newcommand{\bn}{\bf n}
\def\charf {\mbox{{\text 1}\kern-.24em {\text l}}}

\newcommand\Om{\Omega}
\newcommand\ba{\backslash}
\newcommand\BZ{\mathbb Z}
\newcommand\BA{\mathbb A}
\newcommand\BQ{\mathbb Q}

\newcommand\BD{\mathbb D}
\newcommand\BH{\mathbb H}
\newcommand\BR{\mathbb R}
\newcommand\BC{\mathbb C}
\newcommand\lrt{\longrightarrow}
\newcommand\lmt{\longmapsto}
\newcommand\CX{{\Cal X}}
\newcommand\td{\bigtriangledown}
\newcommand\pdx{ {{\partial}\over{\partial x}} }
\newcommand\pdy{ {{\partial}\over{\partial y}} }
\newcommand\pdu{ {{\partial}\over{\partial u}} }
\newcommand\pdv{ {{\partial}\over{\partial v}} }
\newcommand\PZ{ {{\partial}\over {\partial Z}} }
\newcommand\PW{ {{\partial}\over {\partial W}} }
\newcommand\PZB{ {{\partial}\over {\partial{\overline Z}}} }
\newcommand\PWB{ {{\partial}\over {\partial{\overline W}}} }
\newcommand\PX{ {{\partial}\over{\partial X}} }
\newcommand\PY{ {{\partial}\over {\partial Y}} }
\newcommand\PU{ {{\partial}\over{\partial U}} }
\newcommand\PV{ {{\partial}\over{\partial V}} }
\renewcommand\th{\theta}
\renewcommand\l{\lambda}
\renewcommand\k{\kappa}
\newcommand\G{\Gamma}
\newcommand\s{\sigma}
\newcommand\g{\gamma}
\newcommand\PWE{ \frac{\partial}{\partial \overline W}}
\newcommand\fh{\mathfrak h}
\newcommand\fl{\mathfrak l}
\newcommand\wlm{W_{{\mathfrak l},m}}
\newcommand\clm{c_{{\mathfrak l},m}}
\newcommand\rlm{R_{{\mathfrak l},m}}
\newcommand\trlm{{\tilde R}_{{\mathfrak l},m}}
\newcommand\hlm{H_{{\mathfrak l},m}}
\newcommand\slm{s_{{\mathfrak l},m}}
\newcommand\glm{G_{{\mathfrak l},m}}
\newcommand\gtlm{G_{2,{\mathfrak l},m}}
\newcommand\rtlm{R_{2,{\mathfrak l},m}}

\newcommand\rhm{\rho_{\mathcal M}}
\newcommand\hm{H_{\mathcal M}}
\newcommand\CM{\mathcal M}
\newcommand\tgm{\tau^g_{\mathcal M}}
\newcommand\RM{R_{\mathcal M}}
\newcommand\cm{c_{\mathcal M}}
\newcommand\sm{s_{\mathcal M}}
\newcommand\rtm{R_{2,{\mathcal M}}}
\newcommand\pim{\pi_{\mathcal M}}
\newcommand\ztm{Z_{2,\mathcal M}}
\newcommand\hrnm{H_\BR^{(n,m)}}
\newcommand\symm{\textrm{Symm}^2(\BR^m)}
\newcommand\gtm{G_{2,{\mathcal M}}}
\newcommand\symn{\textrm{Symm}^2(\BR^n)}

\newcommand\Jm{J_{\mathcal M}}
\newcommand\dnm{d\nu_{\mathcal M,\Omega}}
\newcommand\uumo{{\mathscr U}_{\mathcal M,\Omega}}
\newcommand\umo{U_{\mathcal M,\Omega}}
\newcommand\hmo{H_{\mathcal M,\Omega}}
\newcommand\pimo{\pi_{\mathcal M,\Omega}}
\newcommand\cmo{c_{\mathcal M,\Omega}}
\newcommand\wm{{\mathscr W}_{\mathcal M}}
\newcommand\wg{\widetilde g}
\newcommand\wgam{\widetilde \gamma}
\newcommand\mwm{{\mathscr W}_{\mathcal M}}

\newcommand\wmo{{\mathscr W}_{\mathcal M,\Omega}}
\newcommand\rmn{\BR^{(m,n)}}

\newcommand\mfm{{\mathscr F}^{(\CM)} }
\newcommand\mfoz{{\mathscr F}^{(\CM)}_{\Om,Z} }
\newcommand\wgm{{\widetilde \gamma} }
\newcommand\Tm{\Theta^{(\CM)} }

\vskip 0.5cm
%%%%%%%%%%%%%%%%%%%%%%%%%%%%%%%%%%%%%%%%%%%%%%%%%%%%%%%%%%%%%%%%%%%%%%%%%%%%%%%%%%%%%%%%%%%%%%%%%%%%%%%%%%%%%%%%%%%%%%%%%%%%%%%%%%%%
%%%%%%%%%%%%%%%%%%%%%%%%%%%%%%%%%%%%%%%%%%%%%%%%%%%%%%%%%%%%%%%%%%%%%%%%%%%%%%%%%%%%%%%%%%%%%%%%%%%%%%%%%%%%%%%%%%%%%%%%%%%%%%%%%%%%
%%%%%%%%%%%%%%%%%%%%%%%%%%%%%%%%%%%%%%%%%%%%%%%%%%%%%%%%%%%%%%%%%%%%%%%%%%%%%%%%%%%%%%%%%%%%%%%%%%%%%%%%%%%%%%%%%%%%%%%%%%%%%%%%%%%%
%%
%%
%%           Section 1       Introduction
%%
%%
%%%%%%%%%%%%%%%%%%%%%%%%%%%%%%%%%%%%%%%%%%%%%%%%%%%%%%%%%%%%%%%%%%%%%%%%%%%%%%%%%%%%%%%%%%%%%%%%%%%%%%%%%%%%%%%%%%%%%%%%%%%%%%%%%%%%
%%%%%%%%%%%%%%%%%%%%%%%%%%%%%%%%%%%%%%%%%%%%%%%%%%%%%%%%%%%%%%%%%%%%%%%%%%%%%%%%%%%%%%%%%%%%%%%%%%%%%%%%%%%%%%%%%%%%%%%%%%%%%%%%%%%%
%%%%%%%%%%%%%%%%%%%%%%%%%%%%%%%%%%%%%%%%%%%%%%%%%%%%%%%%%%%%%%%%%%%%%%%%%%%%%%%%%%%%%%%%%%%%%%%%%%%%%%%%%%%%%%%%%%%%%%%%%%%%%%%%%%%%

\begin{section}{{\large\bf Introduction}}
\setcounter{equation}{0}

\vskip 0.21cm The Weil representation of the symplectic group was
first introduced by A. Weil in his remarkable paper \cite{W} to
reformulate Siegel's analytic theory of quadratic forms\,\cite{Si}
in group theoretical terms. The Weil representation plays a
central role in the study of the transformation behaviors of theta
series and has many applications to the theory of automorphic
forms\,(cf.\,\cite{Ge, KS1, KS2, KS3, LV, Ni, S, Sh}).
The Jacobi group is defined to be the semi-direct product of the
symplectic group and the Heisenberg group.
The Jacobi group is an important object in the framework of
quantum mechanics, geometric quantization and optics
\cite{AAG, B1,B2,B3, GS1, GS2, H, Lu, Per, St, Wo, Yu}.
The squeezed states
in quantum optics represent a physical realization of the coherent states associated
with the Jacobi group \cite{H, Lu, St, Yu}. In this paper, we show that we can construct several types of
the Weil representations of the Jacobi group and present their applications to the theory of
automorphic forms on the Jacobi group and representation theory of
the Jacobi group.

\vskip 0.12cm For a given fixed positive integer $n$, we let
$${\mathbb H}_n=\,\{\,\Om\in \BC^{(n,n)}\,|\ \Om=\,^t\Om,\ \ \ \text{Im}\,\Om>0\,\}$$
be the Siegel upper half plane of degree $n$ and let
$$Sp(n,\BR)=\{ g\in \BR^{(2n,2n)}\ \vert \ ^t\!gJ_ng= J_n\ \}$$
be the symplectic group of degree $n$, where $F^{(k,l)}$ denotes
the set of all $k\times l$ matrices with entries in a commutative
ring $F$ for two positive integers $k$ and $l$, $^t\!M$ denotes
the transposed matrix of a matrix $M,\ \text{Im}\,\Om$ denotes the
imaginary part of $\Om$ and
$$J_n=\begin{pmatrix} 0&I_n \\
                   -I_n&0 \\ \end{pmatrix}.$$
We see that $Sp(n,\BR)$ acts on $\BH_n$ transitively by
\begin{equation*}g\cdot \Om=(A\Om+B)(C\Om+D)^{-1}, \end{equation*}
where $g=\begin{pmatrix} A&B\\ C&D\end{pmatrix}\in Sp(n,\BR)$ and
$\Om\in \BH_n.$

For two positive integers $n$ and $m$, we consider the Heisenberg
group
$$H_{\BR}^{(n,m)}=\{\,(\l,\mu;\k)\,|\ \l,\mu\in \BR^{(m,n)},\ \k\in \BR^{(m,m)},\ \
\k+\mu\,^t\l\ \text{symmetric}\ \}$$ endowed with the following
multiplication law
$$(\l,\mu;\k)\circ (\l',\mu';\k')=(\l+\l',\mu+\mu';\k+\k'+\l\,^t\mu'-
\mu\,^t\l').$$ We let
$$G^J=Sp(n,\BR)\ltimes H_{\BR}^{(n,m)}\quad \ ( \textrm{semi-direct product})$$
be the Jacobi group endowed with the following multiplication law
$$\Big(g,(\lambda,\mu;\kappa)\Big)\cdot\Big(g',(\lambda',\mu';\kappa')\Big) =\,
\Big(gg',(\widetilde{\lambda}+\lambda',\widetilde{\mu}+ \mu';
\kappa+\kappa'+\widetilde{\lambda}\,^t\!\mu'
-\widetilde{\mu}\,^t\!\lambda')\Big)$$ with $g,g'\in Sp(n,\BR),
(\lambda,\mu;\kappa),\,(\lambda',\mu';\kappa') \in
H_{\BR}^{(n,m)}$ and
$(\widetilde{\lambda},\widetilde{\mu})=(\lambda,\mu)g'$.

Then we have the {\it natural action}
of $G^J$ on the Siegel-Jacobi space $\BH_{n,m}:=\BH_n\times
\BC^{(m,n)}$ defined by
\begin{equation}\Big(g,(\lambda,\mu;\kappa)\Big)\cdot (\Om,Z)=\Big(g\!\cdot\! \Om,(Z+\lambda \Om+\mu)
(C\Om+D)^{-1}\Big), \end{equation}

\noindent where $g=\begin{pmatrix} A&B\\
C&D\end{pmatrix} \in Sp(n,\BR),\ (\lambda,\mu; \kappa)\in
H_{\BR}^{(n,m)}$ and $(\Om,Z)\in \BH_{n,m}.$ We refer to
\cite{Y6}-\cite{Y10},\,\cite{Y13}-\cite{Y15} for more details on
materials related to the Siegel-Jacobi space.

\vskip 0.12cm The aim of this article is to introduce three types
of the Weil representations of the Jacobi group $G^J$ and to study
their applications to the theory of automorphic forms and
representation theory. They are slightly different each other. They
are essentially isomorphic. However each has its own advantage in applications to the theory
of automorphic forms and representation theory.

\vskip 0.2cm
This article is
organized as follows. In Section 2, we review the Weil
representation of the symplectic group and the Maslov index briefly. In Section 3, we
define the Weil representation of the Jacobi group $G^J$ using a cocycle class of $G^J$
in $H^2(G^J,T)$.
In Section 4, we define the Schr{\"o}dinger-Weil
representation of the Jacobi group that is used to study the
transformation behaviors of certain theta series with toroidal
variables. The the Schr{\"o}dinger-Weil
representation plays an important role in the construction of Jacobi forms,
the theory of Maass-Jacobi forms and the study of Jacobi's theta sums.
We deal with these applications in detail in Section 7.
In Section 5, we recall the Weil-Satake representation
of the Jacobi group formulated by Satake \cite{Sat} on the Fock model of the
Heisenberg group. In Section 6,
we recall the concept of Jacobi forms of half integral weight to be used in a subsequent section.
We review Siegel modular forms of half integral weight.
In Section 7, we present
the applications of the Schr{\"o}dinger-Weil representation to
constructing of Jacobi forms via covariant maps for the Schr{\"o}dinger-Weil representation,
the study of Maass-Jacobi forms and Jacobi's theta sums. We describe the works of the author \cite{Y17},
A. Piale \cite{P} and J. Marklof \cite{Ma}.
In Section 8, we provides some applications of the
Weil-Satake representation of $G^J$ to the study of
representations of $G^J$ which were obtained by Takase \cite{Ta1, Ta2, Ta4}.
Takase \cite{Ta1} showed that there is a bijective correspondence between the unitary equivalence classes of
unitary representations of a two-fold covering group of the symplectic group and
the unitary equivalence classes of
unitary representations of the Jacobi group.
Using this representation theoretical fact, Takase \cite{Ta4} established a bijective correspondence between the space of cuspidal
Jacobi forms and the space of Siegel cusp forms of half integral weight which is compatible with the action of
Hecke operators.

%In Section 9, we discuss the connection of the Weil
%representations of the Jacobi group to a lifting problem of cusp
%forms.

\vskip 0.2cm \noindent {\bf Notations\,:} \ \ We denote by $\BR$
and $\BC$ the field of real numbers, and the field of complex numbers
respectively. We
denote by $\BR^*_+$ the multiplicative group of positive real
numbers. $\BC^*$ (resp. $\BR^*$) denotes the multiplicative group
of nonzero complex (resp. real) numbers. We denote by $\BZ$ and
$\BZ^+$ the ring of integers and the set of all positive integers
respectively.
$T=\left\{ z\in\BC\,|\ |z|=1\,\right\}$
denotes the multiplicative group of
complex numbers of modulus one. The symbol ``:='' means that the
expression on the right is the definition of that on the left. For
two positive integers $k$ and $l$, $F^{(k,l)}$ denotes the set of
all $k\times l$ matrices with entries in a commutative ring $F$.
For a square matrix $A\in F^{(k,k)}$ of degree $k$, $\sigma(A)$
denotes the trace of $A$. For any $M\in F^{(k,l)},\ ^t\!M$ denotes
the transposed matrix of $M$. $I_n$ denotes the identity matrix of
degree $n$. For a positive integer $m$ we denote by $S(m)$ the set of all
$m\times m$ symmetric real matrices.
We put $i=\sqrt{-1}.$ For $z\in\BC,$ we define
$z^{1/2}=\sqrt{z}$ so that $-\pi / 2 < \ \arg (z^{1/2})\leqq
\pi/2.$ Furthermore we put $z^{\kappa/2}=\big(z^{1/2}\big)^\kappa$ for
every $\kappa\in\BZ.$ For a rational number field $\BQ$, we denote by ${\mathbb A}$
and ${\mathbb A}^*$ the ring of adeles of $\BQ$ and the
multiplicative group of ideles of $\BQ$ respectively.

\end{section}

\vskip 1cm
%%%%%%%%%%%%%%%%%%%%%%%%%%%%%%%%%%%%%%%%%%%%%%%%%%%%%%%%%%%%%%%%%%%%%%%%%%%%%%%%%%%%%%%%%%%%%%%%%%%%%%%%%%%%%%%%%%%%%%%%%%%%%%%%%%%%
%%%%%%%%%%%%%%%%%%%%%%%%%%%%%%%%%%%%%%%%%%%%%%%%%%%%%%%%%%%%%%%%%%%%%%%%%%%%%%%%%%%%%%%%%%%%%%%%%%%%%%%%%%%%%%%%%%%%%%%%%%%%%%%%%%%%
%%%%%%%%%%%%%%%%%%%%%%%%%%%%%%%%%%%%%%%%%%%%%%%%%%%%%%%%%%%%%%%%%%%%%%%%%%%%%%%%%%%%%%%%%%%%%%%%%%%%%%%%%%%%%%%%%%%%%%%%%%%%%%%%%%%%
%%
%%
%%            Section 2  The Weil Representation of the Symplectic Group
%%
%%
%%%%%%%%%%%%%%%%%%%%%%%%%%%%%%%%%%%%%%%%%%%%%%%%%%%%%%%%%%%%%%%%%%%%%%%%%%%%%%%%%%%%%%%%%%%%%%%%%%%%%%%%%%%%%%%%%%%%%%%%%%%%%%%%%%%%
%%%%%%%%%%%%%%%%%%%%%%%%%%%%%%%%%%%%%%%%%%%%%%%%%%%%%%%%%%%%%%%%%%%%%%%%%%%%%%%%%%%%%%%%%%%%%%%%%%%%%%%%%%%%%%%%%%%%%%%%%%%%%%%%%%%%
%%%%%%%%%%%%%%%%%%%%%%%%%%%%%%%%%%%%%%%%%%%%%%%%%%%%%%%%%%%%%%%%%%%%%%%%%%%%%%%%%%%%%%%%%%%%%%%%%%%%%%%%%%%%%%%%%%%%%%%%%%%%%%%%%%%%

\begin{section}{{\large\bf The Weil Representation of the Symplectic Group }}
\setcounter{equation}{0}

\vskip 0.21cm Let $(V,B)$ be a symplectic real
vector space of dimension $2n$ with a non-degenerate alternating
bilinear form $B$. We consider the Lie algebra ${\mathfrak
h}=V+\BR E$ with the Lie bracket satisfying the following properties (2.1) and (2.2)\,:
\begin{equation}
[X,Y]= B(X,Y)E\quad \textrm{for all}\ X,Y \in V;
\end{equation}
\begin{equation}
[Z,E]=0\quad \textrm{for all}\ Z\in {\mathfrak h}.
\end{equation}
Let $H$ be the Heisenberg group with its Lie algebra ${\mathfrak
h}$. Via the exponential map $\exp:{\mathfrak h}\lrt H,\ H$ is
identified with the $(2n+1)$-dimensional vector space with
following multiplication law\,:
\begin{equation*}
\exp (v_1+t_1E)\cdot \exp(v_2+t_2E)=\exp \left(
v_1+v_2+\left(t_1+t_2+{{B(v_1,v_2)}\over 2}\right)E\right),
\end{equation*}
where $v_1,v_2\in V$ and $t_1,t_2\in \BR.$ Let
\begin{equation*}
Sp(B)=\big\{ g\in GL(V)\,|\ B(gx,gy)=B(x,y)\ \quad \textrm{for
all}\ x,y\in V\,\big\}
\end{equation*}
be the symplectic group of $(V,B)$. Then $Sp(B)$ acts on $H$ by
\begin{equation*}
g\cdot \exp (v+tE)= \exp(gv+tE),\quad g\in Sp(B),\ v\in V,\ t\in
\BR.
\end{equation*}
For a fixed nonzero real number $m$, we let $\chi_m:H\lrt T$ be
the function defined by
\begin{equation*}
\chi_m\big(\exp (v+tE)\big)= e^{2\pi i mt},\quad v\in V,\ t\in
\BR.
\end{equation*}
Let ${\mathfrak l}$ be a Lagrangian subspace in $(V,B)$. We put
$L=\exp({\mathfrak l}+\BR E)$. Obviously the restriction of
$\chi_m$ to $L$ is a character of $L$. The induced representation
\begin{equation*}
\wlm= \textrm{Ind}_L^H \chi_m
\end{equation*}
is the so-called $ \textit{Schr{\"o}dinger representation}$ of the
Heisenberg group $H$. The representation $\hlm$ of $\wlm$ is the completion
of the space of continuous functions $\varphi$ on $H$ satisfying the following
properties (2.3) and (2.4)\,:
\begin{equation}
 \varphi(hl)=\chi_m(l)^{-1}\varphi(h),\quad h\in H,\ l\in L\
\end{equation}
and
\begin{equation}
h\mapsto |\varphi(h)|\ \textrm{is square integrable with respect
to an invariant measure on}\ H/L.
\end{equation}
We observe that $\wlm\big(\exp(tE)\big)=e^{2\pi imt}I_{\hlm}$,
where $I_{\hlm}$ denotes the identity operator on $\hlm.$ For
brevity, we put $G=Sp(B).$ For a fixed element $g\in G$, we
consider the representation $\wlm^g$ of $H$ on $\hlm$ defined by
\begin{equation}
\wlm^g(h)=\wlm(g\cdot h),\quad h\in H.
\end{equation}
Since $\wlm(\exp tE)=\wlm^g(\exp tE)$ for all $t\in\BR$, according
to Stone-von Neumann theorem, there exists a unitary operator
$\rlm(g):\hlm\lrt\hlm$ such that $\wlm^g(h)\rlm(g)=\rlm(g)\wlm(h)$
for all $h\in H.$ For convenience, we choose $\rlm({\mathbf
1})=I_{\hlm}$, where ${\mathbf 1}$ denotes the identity element of
$G$. We note that $\rlm(g)$ is determined uniquely up to a scalar
of modulus one. Since $\rlm(g_2)^{-1}\rlm(g_1)^{-1}\rlm(g_1g_2)$
is the unitary operator on $\hlm$ commuting with $\wlm$, according
to Schur's lemma, we have a map $\clm:G\times G\lrt T$ satisfying
the condition
\begin{equation*}
\rlm(g_1g_2)=\clm(g_1,g_2)\rlm(g_1)\rlm(g_2),\quad g_1,g_2\in G.
\end{equation*}
Therefore $\rlm$ is a projective representation of $G$ with
multiplier $\clm$. It is easy to see that the map $\clm$ satisfies
the cocycle condition
\begin{equation*}
\clm(g_1g_2,g_3)\,\clm(g_1,g_2)=\clm(g_1,g_2g_3)\,\clm(g_2,g_3)\quad
\textrm{for all}\ g_1,g_2,g_3\in G.
\end{equation*}
The cocycle $\clm$ produces the central extension $\glm$ of $G$ by
$T$. The group $\glm$ is the set $G\times T$ with the following
group multiplication law\,:
\begin{equation}
(g_1,t_1)\cdot
(g_2,t_2):=\big(g_1g_2,\,t_1t_2\,\clm(g_1,g_2)^{-1}\big),\quad
g_1,g_2\in G,\ t_1,t_2\in T.
\end{equation}
We see that the map $\trlm:\glm\lrt GL(\hlm)$ defined by
\begin{equation*}
\trlm(g,t):=t\,\rlm(g),\quad g\in G,\ t\in\BR
\end{equation*}
is a $ \textit{true}$ representation of $\glm$. \vskip 0.2cm We
now express the cocycle $\clm$ in terms of the Maslov index. Let
$\fl_1,\,\fl_2,\,\fl_3$ be three Lagrangian subspaces of $(V,B)$.
The $ \textit{Maslov index}\ \tau(\fl_1,\fl_2,\fl_3)$ of
$\fl_1,\,\fl_2$ and $\fl_3$ is defined to be the signature of the
quadratic form $Q$ on the $3n$ dimensional vector space
$\fl_1\oplus\fl_2\oplus\fl_3$ given by
\begin{equation*}
Q(x_1+x_2+x_3)=B(x_1,x_2)+B(x_2,x_3)+B(x_3,x_1),\quad x_i\in
\fl_i,\ i=1,2,3.
\end{equation*}
For a sequence $\big\{\fl_1,\fl_2,\cdots,\fl_k\big\}$ of
Lagrangian subspaces $\fl_1,\fl_2,\cdots,\fl_k\,(k\geq 4)$ in
$(V,B)$, we define the $ \textit{Maslov index}\
\tau(\fl_1,\fl_2,\cdots,\fl_k)$ by
\begin{equation*}
\tau(\fl_1,\fl_2,\cdots,\fl_k)=\tau(\fl_1,\fl_2,\fl_3)+\tau(\fl_1,\fl_3,\fl_4)+\cdots+
\tau(\fl_1,\fl_{k-1},\fl_k).
\end{equation*}
For a Lagrangian subspace $\fl$ in $(V,B)$, we put
$\tau_\fl(g_1,g_2)=\tau(\fl,g_1\fl,g_1g_2\fl)$ for $g_1,g_2\in G.$
\begin{lemma}
Let $\fl_1,\fl_2,\cdots,\fl_k$ be Lagrangian subspaces in $(V,B)$
with $k\geq 4$. Then we have
\vskip 0.2cm\noindent (a)
$\tau(\fl_1,\fl_2,\cdots,\fl_k)$ is invariant under the action of
$G$ and its value is unchanged under circular
permutations.\par\noindent (b)\ \ \
$\tau(\fl_1,\fl_2,\fl_3)=-\tau(\fl_2,\fl_1,\fl_3)=-\tau(\fl_1,\fl_3,\fl_2).$
\par\noindent
(c) For any four Lagrangian subspaces $\fl_1,\fl_2,\fl_3,\fl_4$ in
$(V,B)$,
\begin{equation*}
\tau(\fl_1,\fl_2,\fl_3)=\tau(\fl_1,\fl_2,\fl_4)+\tau(\fl_2,\fl_3,\fl_4)+\tau(\fl_3,\fl_1,\fl_4).
\end{equation*}\par\noindent
(d) \ \
$\tau(\fl_1,\fl_2,\cdots,\fl_d)=\tau(\fl_1,\fl_2,\fl)+\tau(\fl_2,\fl_3,\fl)+\cdots+\tau(\fl_{d-1},\fl_d,\fl)
+\tau(\fl_d,\fl_1,\fl)$\ \ for any Lagrangian subspace $\fl$ in
$(V,B)$ and $d\geq 3.$\par\noindent (e)\ \
$\tau(\fl_1,\fl_2,\fl_3,\fl_4)=-\tau(\fl_2,\fl_1,\fl_4,\fl_3).$\par\noindent
(f) For any Lagrangian subspaces
$\fl_1,\fl_2,\fl_3,\fl_1',\fl_2',\fl_3'$ in $(V,B)$, we have
$$\tau(\fl_1',\fl_2',\fl_3')=\tau(\fl_1,\fl_2,\fl_3)+\tau(\fl_1',\fl_2',\fl_2,\fl_1)+\tau(\fl_2',\fl_3',\fl_3,\fl_2)
+\tau(\fl_3',\fl_1',\fl_1,\fl_3).$$ (g) \ \
$\tau_\fl(g_1g_2,g_3)+\tau_\fl(g_1,g_2)=\tau_\fl(g_1,g_2g_3)+\tau_\fl(g_2,g_3)$\
\ for\ all\ $g_1,g_2,g_3\in G$.
\end{lemma}
\noindent {\it Proof.} The proof can be found in \cite{LV}. \hfill
$\square$
\begin{theorem} For a Lagrangian subspace $\fl$ in $(V,B)$ and a
real number $m$, we have
\begin{equation*}
 \clm(g_1,g_2)=e^{-{{i\pi m}\over
4}\,\tau(\fl,g_1\fl,g_1g_2\fl)}\quad for\ all\ g_1,g_2\in G.
\end{equation*}
\end{theorem}
\noindent {\it Proof.} The proof can be found in \cite{LV}. \hfill
$\square$ \vskip 0.2cm An $\textit{oriented vector space}$ of
dimension $n$ is defined to be a pair $(U,e)$, where $U$ is a real
vector space of dimension $n$ and $e$ is an orientation of $U$,
i.e., a connected component of $\bigwedge^n U$--\,$\{0\}.$ For two
oriented vector space $(\fl_1,e_1)$ and $(\fl_2,e_2)$ in a
symplectic vector space $(V,B)$, we define
\begin{equation}
s\big( (\fl_1,e_1),(\fl_2,e_2)\big):=i^{n- \dim
(\fl_1\cap\fl_2)}\,\varepsilon\big( (\fl_1,e_1),(\fl_2,e_2)\big).
\end{equation}
We refer to \cite[pp.\,64--66]{LV} for the precise definition of
$\varepsilon\big( (\fl_1,e_1),(\fl_2,e_2)\big)$. Let $M$ be the
space of all Lagrangian subspaces in $(V,B)$ and $\widetilde M$
the manifold of all oriented Lagrangian subspaces in $(V,B)$. Let
$p:\widetilde M\lrt M$ be the natural projection from $\widetilde
M$ onto $M$. Now we will write $\tilde{\fl}$ for a Lagrangian
oriented subspace $(\fl,e).$
\begin{theorem} Let ${\tilde \fl}_1,\,{\tilde \fl}_2,\,{\tilde
\fl}_3\in \widetilde{M}.$ Then
\begin{equation*}
 e^{-{{i\pi}\over
2}\,\tau(p({\tilde \fl}_1),p({\tilde \fl}_1),p({\tilde
\fl}_1))}=s\big({\tilde \fl}_1,\,{\tilde
\fl}_2\big)\,s\big({\tilde \fl}_2,\,{\tilde
\fl}_3\big)\,s\big({\tilde \fl}_3,\,{\tilde \fl}_1\big).
\end{equation*}
\end{theorem}
\noindent {\it Proof.} The proof can be found in
\cite[pp.\,67--70]{LV}. \hfill $\square$ \vskip 0.2cm Let $\fl$ be
a Lagrangian subspace in $(V,B)$. We choose an orientation $\fl^+$
on $\fl$. Then $G$ acts on oriented Lagrangian subspace in
$(V,B)$. We define
\begin{equation}
\slm(g):=s\big( \fl^+,g\fl^+\big)^m,\quad g\in G.
\end{equation}
The above definition is well defined, i.e., does not depend on the
choice of orientation on $\fl$. Since $\slm(g^{-1})=\slm(g)^{-1},$
according to Theorem 2.1 and Theorem 2.2, we get
\begin{equation}
\clm(g_1,g_2)^2=\, \slm(g_1)^{-1}\,\slm(g_2)^{-1}\,\slm(g_1g_2) \quad
\textrm{for all}\ g_1,g_2\in G.
\end{equation}
Hence we can see that
\begin{equation}
\gtlm:=\big\{\,(g,t)\in \glm\,|\ t^2=\slm(g)^{-1}\,\big\}
\end{equation}
is the subgroup of $\glm$\,(cf.\,Formula (2.6)) that is called
the metaplectic group associated with a pair $(\fl,m)$. We know
that $\gtlm$ is a two-fold covering group of $G$. The restriction
$\rtlm$ of $\trlm$ to $\gtlm$ is a true representation of $\gtlm$
that is called the $ \textit{Weil representation}$ of $G$
associated with a pair $(\fl,m)$. We note that
\begin{equation}
\rtlm(g,t)=t\, \rlm (g)=\slm (g)^{-1/2}\rlm(g)\quad \textrm{for
all}\ (g,t)\in \gtlm.
\end{equation}

We refer to \cite{Ge, KV, LV} for more detail on the Weil representation.

%\begin{eqnarray}
%[X,Y]= B(X,Y)E\quad \textrm{for all}\ X,Y \in V;\\
%[Z,E]=0\quad \textrm{for all}\ Z\in {\mathfrak h}.
%\end{eqnarray}

\end{section}

\vskip 1cm
%%%%%%%%%%%%%%%%%%%%%%%%%%%%%%%%%%%%%%%%%%%%%%%%%%%%%%%%%%%%%%%%%%%%%%%%%%%%%%%%%%%%%%%%%%%%%%%%%%%%%%%%%%%%%%%%%%%%%%%%%%%%%%%%%%%%
%%%%%%%%%%%%%%%%%%%%%%%%%%%%%%%%%%%%%%%%%%%%%%%%%%%%%%%%%%%%%%%%%%%%%%%%%%%%%%%%%%%%%%%%%%%%%%%%%%%%%%%%%%%%%%%%%%%%%%%%%%%%%%%%%%%%
%%%%%%%%%%%%%%%%%%%%%%%%%%%%%%%%%%%%%%%%%%%%%%%%%%%%%%%%%%%%%%%%%%%%%%%%%%%%%%%%%%%%%%%%%%%%%%%%%%%%%%%%%%%%%%%%%%%%%%%%%%%%%%%%%%%%
%%
%%
%%                           Section 3  The Weil Representation of the Jacobi Group $G^J$
%%
%%
%%%%%%%%%%%%%%%%%%%%%%%%%%%%%%%%%%%%%%%%%%%%%%%%%%%%%%%%%%%%%%%%%%%%%%%%%%%%%%%%%%%%%%%%%%%%%%%%%%%%%%%%%%%%%%%%%%%%%%%%%%%%%%%%%%%%
%%%%%%%%%%%%%%%%%%%%%%%%%%%%%%%%%%%%%%%%%%%%%%%%%%%%%%%%%%%%%%%%%%%%%%%%%%%%%%%%%%%%%%%%%%%%%%%%%%%%%%%%%%%%%%%%%%%%%%%%%%%%%%%%%%%%
%%%%%%%%%%%%%%%%%%%%%%%%%%%%%%%%%%%%%%%%%%%%%%%%%%%%%%%%%%%%%%%%%%%%%%%%%%%%%%%%%%%%%%%%%%%%%%%%%%%%%%%%%%%%%%%%%%%%%%%%%%%%%%%%%%%%

\begin{section}{{\large\bf The Weil Representation of the Jacobi Group $G^J$}}
\setcounter{equation}{0}

\vskip 0.21cm Let $V=\BR^{(m,n)}\times \BR^{(m,n)}$ be the symplectic real veactor space with a nondegenerate alternating
bilinear form on $V$ given by
\begin{equation*}
B((\lambda,\mu),(\lambda',\mu')):=\,\sigma( \lambda\,{}^t\!\mu'-\mu\,{}^t\!\lambda'\,),\quad
(\lambda,\mu),(\lambda',\mu')\in\BR^{(m,n)}.
\end{equation*}
We assume that $\CM$ is a positive definite
symmetric real matrix of degree $m$. We denote by $ S(m)$ the set
of all $m\times m$ symmetric real matrices. We let
\begin{equation}
\wm:H_\BR^{(n,m)}\lrt U(\hm)
\end{equation}
be the Schr{\"o}dinger representation with central character
$\wm((0,0;\k))=e^{2\pi i\,\s(\CM \k)}I_{\hm},$\\ $\k\in S(m)$.
Here $\hm$ denotes the representation space of $\wm$.
We note that $\wm$ is
realized on $L^2\big(\BR^{(m,n)}\big)\cong \hm$ by
\begin{equation}
\left( \wm (h)f\right) (x)=\,e^{2\,\pi\,i\,\sigma\left(\CM
(\k\,+\,\mu\,{}^t\!\lambda\,+\,2\,x\,{}^t\mu)\right) }
f(x+\l),\quad x\in \BR^{(m,n)},
\end{equation}

\noindent where $h=(\lambda,\mu;\k)\in \hrnm$ and $f\in
L^2\big(\BR^{(m,n)}\big)$. We refer to \cite{Y1, Y2, Y4, Y5} for more detail about
$\wm$. The Jacobi group $G^J$ acts on $\hrnm$
by conjugation inside $G^J$. Fix an element $\wg\in G^J.$ The
irreducible unitary representation $\mwm^{\wg}$ of $\hrnm$ defined
by
\begin{equation}
\mwm^{\wg} (h):=\mwm(\wg\, h\, \wg^{-1}),\quad h\in \hrnm
\end{equation}
has the property that $\mwm^{\wg}
((0,0;\k))=\mwm((0,0;\k))=e^{2\pi i \s(\CM \k)}\cdot I_{\hm}$ for
all $\kappa\in S(m)$. According to Stone-von Neumann theorem, there
exists a unitary operator $T_\CM (\wg)$ on $\hm$ such that
$T_\CM(\wg)\,\mwm(h)=\mwm^{\wg} (h)\, T_\CM (\wg)$ for all
$h\in\hrnm$. We observe that $T_\CM (\wg)$ is determined uniquely
up to a scalar of modulus one. According to Schur's lemma, we have
a map ${\widetilde c}_\CM:G^J\times G^J\lrt T$ satisfying the
relation
\begin{equation}
T_\CM(\wg_1\wg_2)={\widetilde c}_\CM (\wg_1,\wg_2)\,T_\CM
(\wg_1)T_\CM (\wg_2),\quad \wg_1,\wg_2\in G^J.
\end{equation}
Therefore $T_\CM$ is a projective representation of $G^J$ and
${\widetilde c}_\CM$ defines the cocycle class in $H^2(G^J,T).$
The cocycle ${\widetilde c}_\CM$ satisfies the following properties
\begin{equation}
{\widetilde c}_\CM(h_1,h_2)=1 \quad \textrm{for all}\ h_1,h_2\in H_\BR^{(n,m)};
\end{equation}
\begin{equation}
{\widetilde c}_\CM({\tilde g},e)={\widetilde c}_\CM(e,{\tilde g}) =
{\widetilde c}_\CM(e,e)=1 \quad \textrm{for all}\ {\tilde g}\in G^J\,;
\end{equation}
\begin{equation}
{\widetilde c}_\CM({\tilde g},{\tilde g}^{-1})={\widetilde c}_\CM({\tilde g}^2,{\tilde g}^{-1})\,
{\widetilde c}_\CM( {\tilde g},{\tilde g}) \quad \textrm{for all}\ {\tilde g}\in G^J\,;
\end{equation}
\begin{equation}
T_\CM ( {\tilde g}^{-1})=\, {\widetilde c}_\CM({\tilde g},{\tilde g}^{-1})^{-1}\, T_\CM({\tilde g})^{-1}
\quad \textrm{for all}\ {\tilde g}\in G^J,
\end{equation}
where $e$ is the identity element of $G^J$. The cocycle ${\widetilde c}_\CM$ yields the central extension
$G_\CM^J$ of $G^J$ by $T$. The extension group $G_\CM^J$ is the
set $G^J\times T$ with the following group multiplication law\,:
\begin{equation}
(\wg_1,t_1)\cdot (\wg_2,t_2)=\big(\wg_1\wg_2,t_1t_2\,{\widetilde
c}_\CM (\wg_1,\wg_2)^{-1}\big),\quad \wg_1,\wg_2\in G^J,\
t_1,t_2\in T.
\end{equation}
It is obvious that $\big( (I_{2n},(0,0;0)),1\big)$ is the identity element of $G_\CM^J$ and
$$({\tilde g},t)^{-1}=
\big({\tilde g}^{-1}, t^{-1}\,{\widetilde c}_\CM({\tilde g},{\tilde g}^{-1})\big)$$
if $({\tilde g},t)\in G_\CM^J.$
We see easily that the map ${\widetilde T}_\CM:G_\CM^J\lrt U(\hm)$
defined by
\begin{equation}
{\widetilde T}_\CM (\wg,t):=\,t\,T_\CM (\wg),\quad (\wg,t)\in
G_\CM^J
\end{equation}
is a true representation of $G_\CM^J.$ Here $U(\hm)$ denotes the group of unitary operators of
$\hm$. For the Lagrangian
subspace $\fl=\,\left\{ (0,\mu)\in V\,|\ \mu\in\BR^{(m,n)}\,\right\}$, as (2.8) and (2.9) in Section 2, we can define the
function ${\widetilde s}_\CM:G^J\lrt T$ satisfying the relation
\begin{equation}
{\widetilde c}_\CM(\wg_1,\wg_2)^2=\, {\widetilde
s}_\CM(\wg_1)^{-1}\, {\widetilde s}_\CM (\wg_2)^{-1} \,{\widetilde
s}_\CM(\wg_1 \wg_2)\quad \textrm{for all}\ \ \wg_1,\wg_2\in G^J.
\end{equation}

\noindent Then it is easily seen that
\begin{equation}
G_{\CM,2}^J:=\big\{\,(\wg,t)\in G_\CM^J\,|\ t^2={\widetilde s}_\CM
(\wg)^{-1}\,\big\}
\end{equation}

\noindent is a two-fold covering group of $G^J$. The restriction
${\widetilde \omega}_\CM$ of ${\widetilde T}_\CM$ to $G_{\CM,2}^J$
is called the $\textit{Weil representation}$ of $G^J$ associated
with $\CM$.

\end{section}

\vskip 1cm
%%%%%%%%%%%%%%%%%%%%%%%%%%%%%%%%%%%%%%%%%%%%%%%%%%%%%%%%%%%%%%%%%%%%%%%%%%%%%%%%%%%%%%%%%%%%%%%%%%%%%%%%%%%%%%%%%%%%%%%%%%%%%%%%%%%%
%%%%%%%%%%%%%%%%%%%%%%%%%%%%%%%%%%%%%%%%%%%%%%%%%%%%%%%%%%%%%%%%%%%%%%%%%%%%%%%%%%%%%%%%%%%%%%%%%%%%%%%%%%%%%%%%%%%%%%%%%%%%%%%%%%%%
%%%%%%%%%%%%%%%%%%%%%%%%%%%%%%%%%%%%%%%%%%%%%%%%%%%%%%%%%%%%%%%%%%%%%%%%%%%%%%%%%%%%%%%%%%%%%%%%%%%%%%%%%%%%%%%%%%%%%%%%%%%%%%%%%%%%
%%
%%
%%                    Section 4    The Schr{\"o}dinger-Weil Representation
%%
%%
%%%%%%%%%%%%%%%%%%%%%%%%%%%%%%%%%%%%%%%%%%%%%%%%%%%%%%%%%%%%%%%%%%%%%%%%%%%%%%%%%%%%%%%%%%%%%%%%%%%%%%%%%%%%%%%%%%%%%%%%%%%%%%%%%%%%
%%%%%%%%%%%%%%%%%%%%%%%%%%%%%%%%%%%%%%%%%%%%%%%%%%%%%%%%%%%%%%%%%%%%%%%%%%%%%%%%%%%%%%%%%%%%%%%%%%%%%%%%%%%%%%%%%%%%%%%%%%%%%%%%%%%%
%%%%%%%%%%%%%%%%%%%%%%%%%%%%%%%%%%%%%%%%%%%%%%%%%%%%%%%%%%%%%%%%%%%%%%%%%%%%%%%%%%%%%%%%%%%%%%%%%%%%%%%%%%%%%%%%%%%%%%%%%%%%%%%%%%%%

\begin{section}{{\large\bf The Schr{\"o}dinger-Weil Representation}}
\setcounter{equation}{0}

\vskip 0.21cm Let $\mwm$ be the Schr{\"o}dinger representation of
$\hrnm$ defined by (3.1) in Section 3. The symplectic group
$G=Sp(n,\BR)$ acts on $\hrnm$ by conjugation inside $G^J.$ We fix
an element $g\in G.$ We consider the unitary representation
$\mwm^g$ of $\hrnm$ defined by
\begin{equation}
\mwm^g(h)=\,\mwm (ghg^{-1}), \quad \ \ h\in \hrnm.
\end{equation}

\noindent Since $\mwm^g\big( (0,0;\k)\big)=\,\mwm\big(
(0,0;\k)\big)=\,e^{2\,\pi\,i\,\sigma (\CM\k)} I_{\hm}$ for all
$\k\in S(m),$ according to Stone-von Neumann theorem, $\mwm^g$ is
unitarily equivalent to $\mwm$. Thus there exists a unitary
operator $\RM (g)$ of $\hm$ satisfying the commutation relation
$\RM (g)\,\mwm (h)=\,\mwm^g (h)\,\RM (g)$ for all $h\in\hrnm.$ We
observe that $\RM$ is determined uniquely up to a scalar of
modulus one. According to Schur's lemma, we have a map
$\cm:G\times G\lrt T$ satisfying the relation
\begin{equation}
\RM(g_1g_2)=\cm (g_1,g_2)\,\RM(g_1)\RM (g_2),\quad g_1,g_2\in G.
\end{equation}
Therefore $\RM$ is a projective representation of $G$ and $\cm$
defines the cocycle class in $H^2(G,T).$ The cocycle $\cm$
gives rise to the central extension $G_\CM$ of $G$ by $T$. The
extension group $G_\CM$ is the set $G\times T$ with the following
group multiplication law\,:
\begin{equation}
(g_1,t_1)\cdot (g_2,t_2)=\,\big(g_1g_2,t_1t_2\,c_\CM
(g_1,g_2)^{-1}\big),\quad g_1,g_2\in G,\ t_1,t_2\in T.
\end{equation}
We see that the map ${\widetilde R}_\CM:G_\CM\lrt U(\hm)$ defined
by
\begin{equation}
{\widetilde R}_\CM (g,t)=\,t\,\RM (g),\quad (g,t)\in G_\CM
\end{equation}
is a true representation of $G_\CM.$ For the Lagrangian subspace
$\fl=\,\left\{ (0,\mu)\in V\,|\ \mu\in\BR^{(m,n)}\,\right\}$, as (2.8) and (2.9) in Section 2,
we can define the function
$\sm:G\lrt T$ satisfying the relation
\begin{equation}
c_\CM(g_1,g_2)^2=\sm (g_1)^{-1}\sm (g_2)^{-1} \sm(g_1g_2)\quad
\textrm{for all}\ g_1,g_2\in G.
\end{equation}
Hence we see that
\begin{equation}
\gtm=\big\{\,(g,t)\in G_\CM\,|\ t^2=\sm (g)^{-1}\,\big\}
\end{equation}
is the metaplectic group associated with $\CM\in S(m)$ that is a
two-fold covering group of $G$. The restriction $\rtm$ of
${\widetilde R}_\CM$ to $\gtm$ is the Weil representation of $G$
associated with $\CM\in S(m).$ Now we define the projective
representation $\pim$ of $G^J$ by
\begin{equation}
\pim (hg):=\mwm (h)\,\RM (g),\quad h\in \hrnm,\ g\in G.
\end{equation}

\noindent We observe that any element $\wg$ of $G^J$ can be
expressed in the form $\wg=hg$ with $h\in\hrnm$ and $g\in G.$
Indeed, if $g,g_1\in G$ and $h,h_1\in H_\BR^{(n,m)},$ then we have
\begin{eqnarray*}
\pi_\CM (hgh_1g_1)&=&\,\pi_\CM (hgh_1g^{-1}gg_1)\\
&=&\,\wm (hgh_1g^{-1})\,R_\CM (gg_1)\\
&=&\,c_\CM (g,g_1)\,\wm (h)\,\wm(gh_1g^{-1})\,R_\CM(g) R_\CM(g_1)\\
&=&\,c_\CM (g,g_1)\,\wm (h)\,\wm^g (h_1)\,R_\CM(g) R_\CM(g_1)\\
&=&\,c_\CM (g,g_1)\,\wm (h)\,R_\CM(g)\, \wm(h_1)\,R_\CM(g_1)\\
&=&\,c_\CM (g,g_1)\,\pi_\CM(hg)\,\pi_\CM(h_1g_1).
\end{eqnarray*}
In the second equality, we used the fact that $H_\BR^{(n,m)}$ is a normal
subgroup of $G^J.$ Therefore we get the relation
\begin{equation}
\pi_\CM (hgh_1g_1)=\,c_\CM (g,g_1)\,\pi_\CM(hg)\,\pi_\CM(h_1g_1)
\end{equation}
for all $g,g_1\in G$ and $h,h_1\in H_\BR^{(n,m)}.$ From (4.8) we obtain the relation
\begin{equation}
T_\CM(g)=\,R_\CM (g),\quad {\tilde c}_\CM (g,g')=c_\CM (g,g')\qquad \textrm{for all}
\ g,g'\in G.
\end{equation}

Thus the representation $\pim$ of $G^J$ is naturally extended to the true
representation $\omega^\CM_{\textrm{SW}}$ of $G_{2,\CM}^J:=G_{2,\CM}\ltimes
\hrnm.$ The representation $\omega^\CM_{\textrm{SW}}$ is called
$ \textit{Schr{\"o}dinger}$-{\it Weil representation} of the Jacobi group $G^J$
associated with $\CM\in S(m).$ Indeed we have
\begin{equation}
\omega^\CM_{\textrm{SW}}  (h\!\cdot\!(g,t))=\,t\,\pi_\CM(hg),\quad h\in \hrnm,\
(g,t)\in G_{2,\CM}.
\end{equation}

\noindent
 We recall that the following matrices
\begin{eqnarray*}
t(b):&=&\begin{pmatrix} I_n& b\\
                   0& I_n\end{pmatrix}\ \textrm{with any}\
                   b=\,{}^tb\in \BR^{(n,n)},\\
g(\alpha):&=&\begin{pmatrix} {}^t\alpha & 0\\
                   0& \alpha^{-1}  \end{pmatrix}\ \textrm{with
                   any}\ \alpha\in GL(n,\BR),\\
\s_n:&=&\begin{pmatrix} 0& -I_n\\
                   I_n&\ 0\end{pmatrix}
\end{eqnarray*}
\noindent generate the symplectic group $G=Sp(n,\BR)$
(cf.\,\cite[p.\,326]{F},\,\cite[p.\,210]{Mum}).

\indent The Weil representation $\rtm$ is realized on the Hilbert
space $L^2\big(\BR^{(m,n)}\big)$\,(cf.\,\cite{W},\,\cite{KV})\,:

\begin{equation}
\left( R_\CM (t(b)f\right)(x)=\,e^{{ 2\,\pi\,i\, }\s
(\CM\,x\,b\,{}^tx)}f(x),\quad b=\,{}^tb\in \BR^{(n,n)}\,;
\end{equation}

\begin{equation}
\left( R_\CM (g(\alpha)f\right)(x)=(\det \alpha)^{\frac
m2}f(x\,{}^t\alpha),\quad \alpha\in GL(n,\BR),
\end{equation}

\begin{equation}
\left( R_\CM (\sigma_n)f\right)(x)=\,\left( {1\over
i}\right)^{\frac {mn}2} (\det \CM)^{\frac n2}\int_{\BR^{(m,n)}}
e^{-\,4\,\pi\,i\,\s (\CM\,y\,^t\!x)}f(y)dy.
\end{equation}
According to Formulas (4.11)-(4.13), $\rtm$ is decomposed into two
irreducible representations $\rtm^{\pm}$
\begin{equation}
\rtm=\rtm^+\oplus \rtm^-,
\end{equation}
where $\rtm^+$ and $\rtm^-$ are the even Weil representation and
the odd Weil representation respectively. Obviously the center
${\mathcal Z}_{2,\CM}^J$ of $G_{2,\CM}^J$ is given by
\begin{equation*}
{\mathcal Z}_{2,\CM}^J=\big\{ \big( (I_{2n},1),(0,0;\k)\big)\in
G_{2,\CM}^J\,\big\} \cong S(m).
\end{equation*}
We note that $\omega^\CM_{\textrm{SW}}  |_{G_{2,\CM}}=\rtm$ and
$ \omega^\CM_{\textrm{SW}} (h)=\wm (h)$ for all $h\in \hrnm.$

\end{section}

\vskip 1cm
%%%%%%%%%%%%%%%%%%%%%%%%%%%%%%%%%%%%%%%%%%%%%%%%%%%%%%%%%%%%%%%%%%%%%%%%%%%%%%%%%%%%%%%%%%%%%%%%%%%%%%%%%%%%%%%%%%%%%%%%%%%%%%%%%%%%
%%%%%%%%%%%%%%%%%%%%%%%%%%%%%%%%%%%%%%%%%%%%%%%%%%%%%%%%%%%%%%%%%%%%%%%%%%%%%%%%%%%%%%%%%%%%%%%%%%%%%%%%%%%%%%%%%%%%%%%%%%%%%%%%%%%%
%%%%%%%%%%%%%%%%%%%%%%%%%%%%%%%%%%%%%%%%%%%%%%%%%%%%%%%%%%%%%%%%%%%%%%%%%%%%%%%%%%%%%%%%%%%%%%%%%%%%%%%%%%%%%%%%%%%%%%%%%%%%%%%%%%%%
%%
%%
%%              Section 5     The Weil-Satake Representation
%%
%%%%%%%%%%%%%%%%%%%%%%%%%%%%%%%%%%%%%%%%%%%%%%%%%%%%%%%%%%%%%%%%%%%%%%%%%%%%%%%%%%%%%%%%%%%%%%%%%%%%%%%%%%%%%%%%%%%%%%%%%%%%%%%%%%%%
%%%%%%%%%%%%%%%%%%%%%%%%%%%%%%%%%%%%%%%%%%%%%%%%%%%%%%%%%%%%%%%%%%%%%%%%%%%%%%%%%%%%%%%%%%%%%%%%%%%%%%%%%%%%%%%%%%%%%%%%%%%%%%%%%%%%
%%%%%%%%%%%%%%%%%%%%%%%%%%%%%%%%%%%%%%%%%%%%%%%%%%%%%%%%%%%%%%%%%%%%%%%%%%%%%%%%%%%%%%%%%%%%%%%%%%%%%%%%%%%%%%%%%%%%%%%%%%%%%%%%%%%%

\begin{section}{{\large\bf The Weil-Satake Representation}}
\setcounter{equation}{0} \vskip 0.21cm

\vskip 0.215cm In this section we discuss the realization of the Weil
representation on the Fock model and the Weil-Satake representation due to Satake (cf.\,\cite{Sat}).
We follow the notations in Section 3 and Section 4.
For $g=\begin{pmatrix} A&B\\
C&D\end{pmatrix}\in G$, we set
\begin{equation}
J(g,\Omega)=C\Omega +D,\quad \Omega\in \BH_n.
\end{equation}

\noindent Let $\CM$ be an $m\times m$ symmetric real matrix. We
define the map $\Jm:G^J\times \BH_{n,m}\lrt \BC^*$ by
\begin{equation}
\Jm\big(\wg,(\Omega,Z)\big):= e^{2\pi
i\,\s\big(\CM[Z+\l\Omega+\mu](C\Omega+D)^{-1}C\big)}\, \cdot\,
e^{-2\pi i\,\s\big( \CM(\l
\Omega\,{}^t\l+2\l\,{}^tZ+\k+\mu\,{}^t\l )\big)},
\end{equation}
where $\wg=(g,(\l,\mu;\k))\in G^J$ with $g=\begin{pmatrix} A&B\\
C&D\end{pmatrix}\in G$ and $(\lambda,\mu;\kappa)\in \hrnm.$ Here
$M[N]:=\,{}^tNMN$ is a Siegel's notation for two matrices $M$ and
$N$. The $\Jm$ satisfies the cocycle condition
\begin{equation*}
\Jm(\wg_1\,\wg_2,(\Omega,Z))=\Jm(\wg_1,\wg_2\cdot(\Omega,Z))\,\Jm(\wg_2,(\Omega,Z))
\end{equation*}
for all $\wg_1,\wg_2\in G^J$ and $(\Omega,Z)\in\BH_{n,m}.$ We refer \cite{Sat}
and \cite{Y9} for a construction of $J_\CM$.

\vskip 0.2cm
We
introduce the coordinates $(\Omega,Z)$ on $\BH_{n,m}$ and some
notations.
\begin{eqnarray*}
\Omega\,=&\,X\,+\,iY,\quad\ \ X\,=\,(x_{\mu\nu}),\quad\ \
Y\,=\,(y_{\mu\nu})
\ \ \text{real},\\
Z\,=&U\,+\,iV,\quad\ \ U\,=\,(u_{kl}),\quad\ \ V\,=\,(v_{kl})\ \
\text{real},
%[ dU ]=&\wedge_{k,l} du_{kl}\quad \textrm{and} \quad
%[dV]=\wedge_{k,l} dv_{kl}.
\end{eqnarray*}
\begin{equation*}
[ dX ]=\bigwedge_{\mu\leq\nu} dx_{\mu\nu},\quad [ dY
]=\bigwedge_{\mu\leq\nu} dy_{\mu\nu},
\end{equation*}
\begin{equation*}
[ dU ]=\bigwedge_{k,l} du_{kl},\quad [dV]=\bigwedge_{k,l} dv_{kl}.
\end{equation*}

Now we assume that ${\CM}$ is $ \textit{positive definite}$. We
define the function $\k_\CM:\BH_{n,m}\lrt\BR$ by
\begin{equation}
\k_\CM(\Omega,Z):=e^{-4\pi\,\s(\,{}^tV\CM VY^{-1})}.
\end{equation}
We fix an element $\Omega$ in $\BH_n$. We let $\hmo$ be the
complex Hilbert space consisting of all complex valued holomorphic
functions $f$ on $\BC^{(m,n)}$ such that
\begin{equation*}
\int_{\BC^{(m,n)}}|f(Z)|^2 \,d\nu_{\CM,\Omega}(Z) < \infty,
\end{equation*}
where
\begin{equation*}
d\nu_{\CM,\Omega}(Z)=(\det 2\CM)^n(\det
\textrm{Im}\,\Omega)^{-m}\,\k_\CM(\Omega,Z)\,[dU]\wedge[dV].
\end{equation*}

We define an irreducible unitary representation ${\mathscr
U}_{\CM,\Om}$ of $\hrnm$ on $\hmo$ by
\begin{equation}
\big( {\mathscr U}_{\CM,\Om}
(h)f\big)(Z):=\Jm\big(h^{-1},(\Omega,Z)\big)^{-1}\,f(Z-\l
\Omega-\mu),
\end{equation}
where $h=(\l,\mu;\k)\in\hrnm$, $f\in\hmo$ and $Z\in \BC^{(m,n)}.$
It is known that
for any two elements $\Omega_1$ and $\Omega_2$ of $\BH_n$,
${\mathscr U}_{\CM,\Om_1}$ is equivalent to ${\mathscr U}_{\CM,\Om_1}$ (cf. \cite{Sat}).
Therefore ${\mathscr U}_{\CM,\Om}$ is called the {\it Fock representation} of $H_\BR^{(n,m)}$
associated with $\CM$. Clearly ${\mathscr
U}_{\CM,\Om} ((0,0;\k))=e^{-2\pi i\,\s(\CM\k)}.$ According to
Stone-von Neumann theorem, ${\mathscr U}_{\CM,\Om}$ is equivalent
to ${\mathscr W}_{-\CM}$\,\\
(cf.\,Formula (3.1)). Since the
representation ${\mathscr U}_{\CM,\Om}^g\ (g\in G)$ of $\hrnm$
defined by ${\mathscr U}_{\CM,\Om}^g(h)=\,{\mathscr
U}_{\CM,\Om}(ghg^{-1})$ is equivalent to ${\mathscr U}_{\CM,\Om}$,
there exists a unitary operator $\umo (g)$ of $\hmo$ such that
$\umo(g)\uumo(h)= \uumo^g (h)\umo(g)$ for all $h\in \hrnm.$ Thus
we obtain a projective representation $\umo$ of $G$ on $\hmo$ and
a cocycle ${\widehat c}_{\CM,\Om}:G\times G\lrt T$ satisfying the
condition
\begin{equation*}
\umo(g_1g_2)={\widehat
c}_{\CM,\Om}(g_1,g_2)\,\umo(g_1)\,\umo(g_2),\quad g_1,g_2\in G.
\end{equation*}
Now ${\widehat c}_{\CM,\Om}$ and $\umo(g)$ will be determined
explicitly\,(cf.\,\cite{Sat},\,\cite{Ta1}). In fact,
\begin{equation}
{\widehat c}_{\CM,\Om}(g_1,g_2)=\left( { {\gamma
(g_2^{-1}g_1^{-1}\Omega,g_2^{-1}\Omega)}\over
{\gamma(g_1^{-1}\Omega,\Omega)} }\right)^m,
\end{equation}
where
\begin{equation*}
\gamma(\Omega_1,\Omega_2):=\left(\det\left( { {\Omega_1-{\overline
\Omega}_2}\over {2i} }\right)\right)^{-{\frac 12}}\,\left( \det
\rm{Im}\,\Omega_1\right)^{\frac 14}\,\left( \det
\rm{Im}\,\Omega_2\right)^{\frac 14},\quad
\Omega_1,\Omega_2\in\BH_n.
\end{equation*}
We define the projective representation $\tau_{\CM,\Omega}$ of
$G^J$ by
\begin{equation}
\tau_{\CM,\Om}(hg):=\,{\mathscr U}_{\CM,\Om}(h)\,\umo(g)\quad
\textrm{for all}\ h\in\hrnm,\ g\in G.
\end{equation}
Then $\tau_{\CM,\Om}$ satisfies the following relation
\begin{equation}
\tau_{\CM,\Om}(\wg_1\,\wg_2)=\,{\widehat c}_{\CM,\Om}
(g_1,g_2)\,\tau_{\CM,\Om}(\wg_1)\,\tau_{\CM,\Om}(\wg_2)
\end{equation}
for all $\wg_1=(g_1,h_1),\,\wg_2=(g_2,h_2)\in G^J$ with
$g_1,g_2\in
G$ and $h_1,h_2\in \hrnm.$  \\
\indent We put
\begin{equation}
\beta_\Omega(g_1,g_2):= \,{\widehat c}_{\CM,\Om}
\big(g_1,g_2\big)^{{\frac 1m}},\quad g_1,g_2\in G.
\end{equation}
Then $\beta_\Omega$ satisfies the cocycle condition and the
following relation
\begin{equation*}
\beta_\Omega(g_1,g_2)^2= {\widehat s}_\Omega(g_1)^{-1}\,{\widehat
s}_\Omega(g_2)^{-1}\,{\widehat s}_\Omega(g_1g_2),\quad g_1,g_2\in
G,
\end{equation*}
where
\begin{equation*}
{\widehat s}_\Omega (g)=\,|\det J(g^{-1},\Omega)|^{-1}\, \left( \det
J(g^{-1},\Omega)\right),\quad g\in G.
\end{equation*}
The cocycle class $[\beta_\Omega]$ in $H^2(G,T)$ defines the
central extension $G_\Omega=G\times T$ of $G$ by $T$ with the
following multiplication law
\begin{equation*}
(g_1,t_1)\cdot (g_2,t_2)=\big(g_1g_2, t_1t_2\,\beta_\Omega
(g_1,g_2)^{-1}\big).
\end{equation*}
We obtain a normal closed subgroup $G_{2,\Omega}$ of $G_\Omega$
given by
\begin{equation}
G_{2,\Omega}=\left\{\,(g,t)\in G_\Omega\,|\ t^2=\,{\widehat
s}_\Omega(g)^{-1}\,\right\}.
\end{equation}
We can show that $G_{2,\Omega}$ is a two-fold covering group of
$G$. We set for any $g\in G$ and $\Omega_1,\Omega_2\in \BH_n$,
\begin{equation}
\varepsilon(g;\Omega_1,\Omega_2):={
{\gamma(g\!\cdot\!\Omega_1,g\!\cdot\!\Omega_2)}\over
{\gamma(\Omega_1,\Omega_2)} }.
\end{equation}
We can see that for any element $g\in G$ and $\Omega\in\BH_n$, the
topological group $G_{2,\Omega}$ is isomorphic to
$G_{2,g\cdot\Omega}$ via the correspondence
\begin{equation*}
(g_0,t_0)\mapsto \big(
g_0,t_0\,\varepsilon(g_0^{-1};g\!\cdot\!\Omega,\Omega)\big),\quad
(g_0,t_0)\in G_{2,\Omega}.
\end{equation*}
Therefore it is enough to consider only the case $\Omega=iI_n.$ We set
$G_2:=G_{2,\,iI_n}$. We let
\begin{equation*}
G^J_2:=G_2\ltimes \hrnm
\end{equation*}
be the two-fold covering group of $G^J$ endowed with the
multiplication law
\begin{eqnarray*}
&\ \Big((g,t),(\lambda,\mu;\kappa)\Big)\cdot\Big((g',t'),(\lambda',\mu';\kappa')\Big)\\
&=\, \Big(
(g,t)\!\cdot\!(g',t'),(\tilde{\lambda}+\lambda',\tilde{\mu}+ \mu';
\kappa+\kappa'+\tilde{\lambda}\,^t\!\mu'
-\tilde{\mu}\,^t\!\lambda')\Big)
\end{eqnarray*}
with $(g,t),(g',t')\in G_2,\,
(\lambda,\mu;\kappa),(\lambda',\mu';\kappa') \in H_{\BR}^{(n,m)}$
and $(\tilde{\lambda},\tilde{\mu})=(\lambda,\mu)g'$.

\vskip 0.352cm
We note that any element $\widetilde\sigma$ of $G_2^J$ can be written in the form
$\widetilde\sigma=\,h(g,t)$ with $h\in H_\BR^{(n,m)}$ and $(g,t)\in G_2.$
We define a
unitary representation ${\widehat\omega}_{\CM}:={\widehat\omega}_{\CM,iI_n}$ of $G_2^J$ by
\begin{equation}
{\widehat\omega}_{\CM}
(h(g,t)):=\,t^m\,\tau_{\CM,\,iI_n}(hg),\qquad  h\in\hrnm,\ (g,t)\in G_2.
\end{equation}

In fact, if $h,h_1\in H_\BR^{(n,m)}$ and $(g,t),(g_1,t_1)\in G_2,$ then we obtain
\newcommand{\whoc}{{\widehat\omega}_{\CM}}
\begin{eqnarray*}
& & {\widehat\omega}_{\CM} \big(h(g,t)h_1(g_1,t_1) \big)\\
&=&\,\whoc \big(h(g,t)h_1(g,t)^{-1}(g,t)(g_1,t_1) \big)\\
&=&\, \whoc \big(h(g,t)h_1(g,t)^{-1} (gg_1,tt_1\beta_{iI_n}(g,g_1)^{-1}) \big)\\
&=&\, (tt_1)^m \beta_{iI_n}(g,g_1)^{-m} \tau_{\CM,iI_n}
\big( h(g,t)h_1(g,t)^{-1}gg_1 \big)\\
&=&\, (tt_1)^m \beta_{iI_n}(g,g_1)^{-m}
{\mathscr U}_{\CM,iI_n}\big( h(g,t)h_1(g,t)^{-1} \big)
U_{\CM,iI_n}(gg_1) \\
%&=&\, (tt_1)^m {\mathscr U}_{\CM,iI_n}\big( h(g,t)h_1(g,t)^{-1} \big)
%U_{\CM,iI_n}(g) U_{\CM,iI_n}(g_1) \\
&=&\, (tt_1)^m {\mathscr U}_{\CM,iI_n}(h) {\mathscr U}_{\CM,iI_n}^g (h_1)
\,U_{\CM,iI_n}(g) \,U_{\CM,iI_n}(g_1)\\
&=&\, t^m t_1^m \,{\mathscr U}_{\CM,iI_n}(h)\,U_{\CM,iI_n}(g){\mathscr U}_{\CM,iI_n} (h_1)\,
U_{\CM,iI_n}(g_1)\\
&=&\, t^mt_1^m \,\tau_{\CM,\,iI_n}(hg)\,\tau_{\CM,\,iI_n}(h_1g_1)\\
&=&\,{\widehat\omega}_{\CM} \big(h(g,t)\big)\,{\widehat\omega}_{\CM} \big(h_1(g_1,t_1)\big).
\end{eqnarray*}

\vskip 0.21cm\noindent
$\widehat{\omega}_{\CM}$ is called the {\it Weil-Satake
representation} of $G^J$ associated with $\CM$. In Section 8, we discuss some applications of
the Weil-Satake representation ${\widehat\omega}_{\CM}$  to the study of unitary representations of
$G^J$.
%The restriction of
%$\widehat{\omega}_{\CM}$ to $G_2$ is equivalent to the
%$m$-fold tensor product of the usual Weil representation of
%$G$\,(cf.\,Section 2).

%\begin{equation}
%A
%\end{equation}

\vskip 0.21cm

\end{section}

\vskip 1cm
%%%%%%%%%%%%%%%%%%%%%%%%%%%%%%%%%%%%%%%%%%%%%%%%%%%%%%%%%%%%%%%%%%%%%%%%%%%%%%%%%%%%%%%%%%%%%%%%%%%%%%%%%%%%%%%%%%%%%%%%%%%%%%%%%%%%
%%%%%%%%%%%%%%%%%%%%%%%%%%%%%%%%%%%%%%%%%%%%%%%%%%%%%%%%%%%%%%%%%%%%%%%%%%%%%%%%%%%%%%%%%%%%%%%%%%%%%%%%%%%%%%%%%%%%%%%%%%%%%%%%%%%%
%%%%%%%%%%%%%%%%%%%%%%%%%%%%%%%%%%%%%%%%%%%%%%%%%%%%%%%%%%%%%%%%%%%%%%%%%%%%%%%%%%%%%%%%%%%%%%%%%%%%%%%%%%%%%%%%%%%%%%%%%%%%%%%%%%%%
%%
%%
%%               Section 6  Jacobi Forms
%%
%%%%%%%%%%%%%%%%%%%%%%%%%%%%%%%%%%%%%%%%%%%%%%%%%%%%%%%%%%%%%%%%%%%%%%%%%%%%%%%%%%%%%%%%%%%%%%%%%%%%%%%%%%%%%%%%%%%%%%%%%%%%%%%%%%%%
%%%%%%%%%%%%%%%%%%%%%%%%%%%%%%%%%%%%%%%%%%%%%%%%%%%%%%%%%%%%%%%%%%%%%%%%%%%%%%%%%%%%%%%%%%%%%%%%%%%%%%%%%%%%%%%%%%%%%%%%%%%%%%%%%%%%
%%%%%%%%%%%%%%%%%%%%%%%%%%%%%%%%%%%%%%%%%%%%%%%%%%%%%%%%%%%%%%%%%%%%%%%%%%%%%%%%%%%%%%%%%%%%%%%%%%%%%%%%%%%%%%%%%%%%%%%%%%%%%%%%%%%%

\begin{section}{{\large\bf Jacobi Forms}}
\setcounter{equation}{0} \vskip 0.21cm
%\vspace{0.1in}\\
%\indent

\vskip 0.2cm  Let $\rho$ be a rational representation of
$GL(n,\mathbb{C})$ on a finite dimensional complex vector space
$V_{\rho}.$ Let ${\mathcal M}\in \mathbb R^{(m,m)}$ be a symmetric
half-integral semi-positive definite matrix of degree $m$. Let
$C^{\infty}(\BH_{n,m},V_{\rho})$ be the algebra of all
$C^{\infty}$ functions on $\BH_{n,m}$ with values in $V_{\rho}.$
For $f\in C^{\infty}(\BH_{n,m}, V_{\rho}),$ we define
\begin{align}
  & (f|_{\rho,{\mathcal M}}[(g,(\lambda,\mu;\kappa))])(\Om,Z) \notag \\
:= \,& e^{-2\,\pi\, i\,\sigma\left( {\mathcal M}(Z+\lambda
\Om+\mu)(C\Om+D)^{-1}C\,{}^t(Z+\lambda\,\Om\,+\,\mu)\right) }
\times e^{2\,\pi\, i\,\sigma\left( {\mathcal M}(\lambda\,
\Om\,{}^t\!\lambda\,+\,2\,\lambda\,{}^t\!Z+\,\kappa+
\mu\,{}^t\!\lambda) \right)} \\
&\times\rho(C\Om+D)^{-1}f(g\!\cdot\!\Om,(Z+\lambda
\Om+\mu)(C\Om+D)^{-1}),\notag
\end{align}
where $g=\left(\begin{matrix} A&B\\ C&D\end{matrix}\right)\in
Sp(n,\mathbb R),\ (\lambda,\mu;\kappa)\in H_{\mathbb R}^{(n,m)}$
and $(\Om,Z)\in \BH_{n,m}.$
\vspace{0.1in}\\
%\noindent {\bf Definition\ 6.1.}\quad 
\begin{definition}
Let $\rho$ and $\mathcal M$
be as above. Let
$$H_{\mathbb Z}^{(n,m)}:= \{ (\lambda,\mu;\kappa)\in H_{\mathbb R}^{(n,m)}\, \vert
\, \lambda,\mu\in \mathbb Z^{(m,n)},\ \kappa\in \mathbb
Z^{(m,m)}\,\ \}.$$ A {\it Jacobi\ form} of index $\mathcal M$ with
respect to $\rho$ on a subgroup $\Gamma$ of $\Gamma_n$ of finite
index is a holomorphic function $f\in
C^{\infty}(\BH_{n,m},V_{\rho})$ satisfying the following
conditions (A) and (B):

\smallskip

\noindent (A) \,\ $f|_{\rho,{\mathcal M}}[\tilde{\gamma}] = f$ for
all $\tilde{\gamma}\in {\widetilde\Gamma}:= \Gamma \ltimes
H_{\mathbb Z}^{(n,m)}$.

\smallskip

\noindent (B) \,\ For each $M\in\Gamma_n$, $f|_{\rho,\CM}[M]$ has a
Fourier expansion of the following form :
$$\big(f|_{\rho,\CM}[M]\big)(\Om,Z) = \sum\limits_{T=\,{}^tT\ge0\atop \text {half-integral}}
\sum\limits_{R\in \mathbb Z^{(n,m)}} c(T,R)\cdot e^{{ {2\pi
i}\over {\lambda_\G}}\,\sigma(T\Om)}\cdot e^{2\pi i\sigma(RZ)}$$

with a suitable $\lambda_\G\in\BZ$ and
$c(T,R)\ne 0$ only if $\left(\begin{matrix} { 1\over {\lambda_\G}}T & \frac 12R\\
\frac 12\,^t\!R&{\mathcal M}\end{matrix}\right) \geqq 0$.
\end{definition}

\medskip

\indent If $n\geq 2,$ the condition (B) is superfluous by K{\"
o}cher principle\,(\,cf.\,\cite{Zi} Lemma 1.6). We denote by
$J_{\rho,\mathcal M}(\Gamma)$ the vector space of all Jacobi forms
of index $\mathcal{M}$ with respect to $\rho$ on $\Gamma$.
Ziegler\,(\,cf.\,\cite{Zi} Theorem 1.8 or \cite{EZ} Theorem 1.1\,)
proves that the vector space $J_{\rho,\mathcal {M}}(\Gamma)$ is
finite dimensional. In the special case $\rho(A)=(\det(A))^k$ with
$A\in GL(n,\BC)$ and a fixed $k\in\BZ$, we write $J_{k,\CM}(\G)$
instead of $J_{\rho,\CM}(\G)$ and call $k$ the {\it weight} of the
corresponding Jacobi forms. For more results on Jacobi forms with
$n>1$ and $m>1$, we refer to \cite{Y6}-\cite{Y10} and \cite{Zi}. Jacobi forms play an
important role in elliptic cusp forms to Siegel cusp forms of degree $2n$ (cf.\,\cite{I}).
\vskip 0.2cm
%\vspace{0.1in}\\
% \noindent {\bf Definition\ 6.2.}\quad 
\begin{definition} 
A Jacobi form $f\in J_{\rho,\mathcal {M}}(\Gamma)$ is said to be
a {\it cusp}\,(\,or {\it cuspidal}\,) form if $\begin{pmatrix} {
1\over {\lambda_\G}}T & {\frac 12}R\\ {\frac 12}\,^t\!R & \mathcal
{M}\end{pmatrix} > 0$ for any $T,\,R$ with $c(T,R)\ne 0.$ A Jacobi
form $f\in J_{\rho,\mathcal{M}}(\Gamma)$ is said to be {\it
singular} if it admits a Fourier expansion such that a Fourier
coefficient $c(T,R)$ vanishes unless $\text{det}\begin{pmatrix} {
1\over {\lambda_\G}}T &{\frac 12}R\\ {\frac 12}\,^t\!R & \mathcal
{M}\end{pmatrix}=0.$
\end{definition}

\vskip 0.2cm Singular Jacobi forms were characterized by a certain differential operator and the weight by the
author \cite{Y8}.

\vskip 0.53cm
Without loss of generality we may assume that $\rho$ is
irreducible. Then we choose a hermitian inner product $\langle\ ,\
\rangle$ on $V_{\rho}$ that is preserved under the unitary group
$U(n)\subset GL(n,\BC).$ For two Jacobi forms $f_1$ and $f_2$ in
$J_{\rho,\CM}(\Gamma)$, we define the Petersson inner product
formally by
\begin{equation}
\langle f_1,f_2 \rangle:=\int_{\Gamma_{n,m}\backslash \BH_{n,m}}\langle\,
\rho(Y^{\frac 12})f_1(\Omega,Z),\rho(Y^{\frac
12})f_2(\Omega,Z)\rangle\,\k_\CM(\Omega,Z)\,dv,
\end{equation}
where
\begin{equation}
dv=(\det Y)^{-(n+m+1)}[dX]\wedge [dY]\wedge [dU]\wedge [dV]
\end{equation}
is a $G^J$-invariant volume element on $\BH_{n,m}.$
See (5.3) for the definition of $\kappa_\CM (\Om,Z).$
A Jacobi form
$f$ in $J_{\rho,\CM}(\Gamma)$ is said to be square integrable if
$\langle f,f\rangle < \infty.$ We note that cusp Jacobi forms are
square integrable and that $\langle f_1,f_2\rangle$ is finite if
one of $f_1$ and $f_2$ is a cusp Jacobi form (cf.\,\cite{Zi},\
p.\,203).

\smallskip We define the map $J_{\rho,\CM}:G^J\times\BH_{n,m}\lrt
GL(V_\rho)$ by
\begin{equation*}
J_{\rho,\CM}(\wg,(\Om,Z))=J_\CM(\wg,(\Om,Z))\,\rho(J(g,\Omega))\quad\quad
( \textrm{cf.}\,(5.1) \ \textrm{and}\ (5.2)),
\end{equation*}
where $\wg=(g,h)\in G^J$ with $g\in G$ and $h\in\hrnm.$ For a
function $f$ on $\BH_n$ with values in $V_\rho$, we can lift $f$
to a function $\Phi_f$ on $G^J$\,:
\begin{eqnarray*}
\Phi_f(\s):&=& (f|_{\rho,\CM}[\s])(iI_n,0)\\
&=& J_{\rho,\CM}(\s,(iI_n,0))^{-1}f(\s\!\cdot\! (iI_n,0)),\quad
\s\in G^J.
\end{eqnarray*}
A characterization of $\Phi_f$ for a cusp Jacobi form $f$ in
$J_{\rho,\CM}(\Gamma)$ was given by Takase
\cite[pp.\,162--164]{Ta1} and the author
\cite[pp.\,252--254]{Y11}.

\vskip 0.53cm
We allow a weight $k$ to be half-integral. For brevity, we set $G=Sp(n,\BR).$ For any $g\in G$
and $\Omega,\Omega'\in \BH_n$, we note that
\begin{eqnarray}
\varepsilon (g;\Om',\Om)&=&\,\textrm{det} ^{-{\frac 12}}\left( { {g\cdot\Om'-{\overline{g\cdot\Om}} }\over {2\,i}}
\right)\,\textrm{det} ^{{\frac 12}}\left( { {\Om'-{\overline{\Om}} }\over {2\,i}}
\right)\\
& & \ \ \times \,|\det J(g,\Om')|^{-1/2} |\det J(g,\Om)|^{-1/2}.\qquad \quad \quad (\textrm{cf}.\,(5.10))\nonumber
\end{eqnarray}
Here $J(g,\Om)=C\Om+D$ for $g=\begin{pmatrix} A & B \\ C & D \end{pmatrix}
\in G$ (cf. (5.1)).

\vskip 0.2cm
Let $\mathscr S=\,\left\{ S\in\BC^{(n,n)}\,|\ S=\,{}^tS,\ \textrm{Re}(S)>0\,\right\}$ be a connected
simply connected complex manifold. Then there is a uniquely determined holomorphic function
$ \textrm{det}^{\frac 12}$ on $\mathscr S$ such that
\begin{eqnarray}
\big( \textrm{det}^{1/2} S\big)^2&=&\,\det S\qquad \textrm{for all} \ S\in {\mathscr S}.\\
\textrm{det}^{1/2} S &=&\,(\det S)^{1/2} \qquad \textrm{for all}\ S\in {\mathscr S}\cap \BR^{(n,n)}.
\end{eqnarray}

For each integer $k\in\BZ$ and $S\in {\mathscr S}$, we put
$$ \textrm{det}^{k/2}S =\,\big( \textrm{det}^{1/2} S\big)^k.$$
For each $\Om\in \BH_n$, we define the function $\beta_\Om:G\times G\lrt T$ by
\begin{equation}
\beta_\Om (g_1,g_2)=\,\epsilon (g_1;\Om,g_2(\Om)),\qquad g_1,g_2\in G.
\end{equation}
Then $\beta_\Om$ satisfies the cocycle condition and the cohomology class of $\beta_\Om$
of order two\,;
\begin{equation}
\beta_\Om (g_1,g_2)^2=\,\alpha_\Om (g_2)\, \alpha_\Om(g_1g_2)^{-1} \alpha_\Om(g_1),
\end{equation}
where
\begin{equation}
\alpha_\Om (g)=\,{ {\det J(g,\Om)}\over {|\det J(g,\Om)|} },\qquad g\in G,\ \Om\in \BH_n.
\end{equation}
For any $\Om\in \BH_n$, we let
$$G_\Om=\,\left\{ (g,\epsilon)\in G \times T\,|\ \epsilon^2=\,\alpha_\Om(g)^{-1}\,\right\}$$
be the two-fold covering group with multiplication law
\begin{equation*}
(g_1,\epsilon_1) (g_2,\epsilon_2)=\,\big(g_1g_2, \epsilon_1\,\epsilon_2 \beta_\Om (g_1,g_2)\big).
\end{equation*}
The covering group $G_\Om$ depends on the choice of $\Om\in\BH_n$, i.e., the choice of a
maximal compact subgroup of $G$. However for any two elements $\Om_1,\Om_2\in\BH_n,\ G_{\Om_1}$
is isomorphic to $G_{\Om_2}$ (cf.\ \cite{Ta3}). We put $G_*:=\,G_{iI_n}.$

\vskip 0.2cm
We define the automorphic factor $J_{1/2}:G_*\times \BH_n\lrt \BC^*$ by
\begin{equation}
J_{1/2}(g_\epsilon,\Om):=\,\epsilon^{-1}\,\varepsilon (g;\Om,iI_n)\,|\det J(g,\Om)|^{1/2},
\end{equation}
where $g_\epsilon=(g,\epsilon)\in G_\Om$ with $g\in G$ and $\Om\in\BH_n.$ It is easily checked that
\begin{equation}
J_{1/2}(g_*h_*,\Om)=\,J_{1/2}(g_*,h\cdot\Om)J_{1/2}(h_*,\Om)
\end{equation}
for all
$g_*=(g,\epsilon),\,h_*=(h,\eta)\in G_*$ and $\Om\in \BH_n.$
and
\begin{equation}
J_{1/2}(g_*,\Om)^2=\,\det (C\Om+D)
\end{equation}
for all $g_*=(g,\epsilon)\in G$ with
$g=\begin{pmatrix} A & B \\ C & D \end{pmatrix}
\in G.$

\vskip 0.2cm Let $\pi_*:G_*\lrt G$ be the projection defined by $\pi_*(g,\epsilon)=g.$
Let $\G$ be a subgroup of the Siegel modular group $\G_n$ of finite index. Let
$\G_*=\pi_*^{-1}(\G)\subset G_*.$ Let $\chi$ be a finite order unitary character of $\G_*.$
Let $k\in \BZ^+$ be a positive integer.
We say that a holomorphic function $\phi:\BH_n\lrt \BC^*$ is a Siegel modular form of
a half-integral weight $k/2$ with level $\G$ if it satisfies the condition
\begin{equation}
\phi(\gamma_*\cdot\Om)=\,\chi(\g_*)\,J_{1/2}(\g_*,\Om)^k \phi(\Om)
\end{equation}
for all $\gamma_*\in\G_*$ and $\Om\in\BH_n.$ We denote by $M_{k/2}(\G,\chi)$ be the vector space of
all Siegel modular forms of weight $k/2$ with level $\G.$ Let $S_{k/2}(\G,\chi)$ be the subspace
of $M_{k/2}(\G,\chi)$ consisting of $\phi\in M_{k/2}(\G,\chi)$ such that
\begin{equation*}
|\phi(\Om)|\,\det( \textrm{Im}\,\Om)^{k/4}\ \ \textrm{is bounded on}\ \BH_n.
\end{equation*}
An element of $S_{k/2}(\G,\chi)$ is called a Siegel cusp form of weight $k/2.$

\vskip 0.3cm
\noindent 
%{\bf Definition\ 6.3.} 
\begin{definition}
Let $\G\subset \G_n$ be a
subgroup of finite index. We put $\G_*=\pi_*^{-1}(\G)$ and
$$ {\widetilde\G}_*=\,\G_*\ltimes H_\BZ^{(n,m)}.$$
A holomorphic function $f:\BH_{n,m}\lrt
\BC$ is said to be a Jacobi form of a weight $k/2\in {\frac 12}\BZ$ (k\,:\,odd)
with level $\G$ and index $\CM$ for the character $\chi$ of $\Gamma_*$ of if it satisfies the following
transformation formula
\begin{equation}
f({\widetilde\g}_*\cdot
(\Om,Z))=\,\chi(\gamma_*)\,J_{k,\CM}({\widetilde\g}_*,(\Om,Z))f(\Om,Z)\quad
\textrm{for\ all}\ {\widetilde\g}_*\in {\widetilde\G}_*
\end{equation}

\noindent where $J_{k,\CM}:{\widetilde\G}_*\times\BH_{n,m}\lrt\BC$
is an automorphic factor defined by

\begin{eqnarray}
J_{k,\CM}\big( {\widetilde \g}_*,(\Om,Z)\big):=& e^{2\,  \pi\,
i\,\sigma\big({\mathcal{M}}(Z+\lambda
\Om+\mu)(C\Om+D)^{-1}C\,{}^t(Z+\lambda \Om+\mu) \big)} \hskip 2cm
\\ & \times\,\, e^{-2\pi i\sigma\left( {\mathcal{M}}(\lambda
\Om\,{}^t\!\lambda+2\lambda\,{}^t\!Z\,+\,\kappa\,+\,
\mu\,{}^t\!\lambda)\right) } J_{1/2}(\g_*,\Om)^k,\nonumber\hskip 1cm
\end{eqnarray}
where ${\widetilde \g}_*=(\g_*,(\lambda,\mu;\kappa))\in {\widetilde\G}_*$ with $\g=\begin{pmatrix} A&B\\
C&D\end{pmatrix}\in \G,\ \g_*=(\g,\epsilon),\ (\lambda,\mu,\kappa)\in H_{\BZ}^{(n,m)}$
and $(\Om,Z)\in \BH_{n,m}.$
\end{definition}

\end{section}

\vskip 1cm
%%%%%%%%%%%%%%%%%%%%%%%%%%%%%%%%%%%%%%%%%%%%%%%%%%%%%%%%%%%%%%%%%%%%%%%%%%%%%%%%%%%%%%%%%%%%%%%%%%%%%%%%%%%%%%%%%%%%%%%%%%%%%%%%%%%%
%%%%%%%%%%%%%%%%%%%%%%%%%%%%%%%%%%%%%%%%%%%%%%%%%%%%%%%%%%%%%%%%%%%%%%%%%%%%%%%%%%%%%%%%%%%%%%%%%%%%%%%%%%%%%%%%%%%%%%%%%%%%%%%%%%%%
%%%%%%%%%%%%%%%%%%%%%%%%%%%%%%%%%%%%%%%%%%%%%%%%%%%%%%%%%%%%%%%%%%%%%%%%%%%%%%%%%%%%%%%%%%%%%%%%%%%%%%%%%%%%%%%%%%%%%%%%%%%%%%%%%%%%
%%
%%
%%   Section 7   Applications of the Schr{\"o}dinger-Weil Representation
%%
%%%%%%%%%%%%%%%%%%%%%%%%%%%%%%%%%%%%%%%%%%%%%%%%%%%%%%%%%%%%%%%%%%%%%%%%%%%%%%%%%%%%%%%%%%%%%%%%%%%%%%%%%%%%%%%%%%%%%%%%%%%%%%%%%%%%
%%%%%%%%%%%%%%%%%%%%%%%%%%%%%%%%%%%%%%%%%%%%%%%%%%%%%%%%%%%%%%%%%%%%%%%%%%%%%%%%%%%%%%%%%%%%%%%%%%%%%%%%%%%%%%%%%%%%%%%%%%%%%%%%%%%%
%%%%%%%%%%%%%%%%%%%%%%%%%%%%%%%%%%%%%%%%%%%%%%%%%%%%%%%%%%%%%%%%%%%%%%%%%%%%%%%%%%%%%%%%%%%%%%%%%%%%%%%%%%%%%%%%%%%%%%%%%%%%%%%%%%%%

\begin{section}{{\large\bf Applications of the Schr{\"o}dinger-Weil Representation}}
\setcounter{equation}{0}

\vskip 0.5cm \noindent {\bf 7.1. Construction of Jacobi Forms} \vskip 0.3cm

\vskip 0.2cm We assume that $\CM$ is a positive definite symmetric
integral matrix of degree $m$. Let $\omega_\CM$ be the
Schr{\"o}dinger-Weil representation of $G^J$ constructed in
Section 4. We recall that $\omega_\CM$ is realized on the Hilbert
space $L^2\big( \BR^{(m,n)}\big)$ by Formulas (4.9)-(4.11). We
define the mapping ${\mathscr F}^{(\CM)}:\BH_{n,m}\lrt
L^2\big(\rmn\big)$ by

\begin{equation}
\mfm (\Om,Z)(x)=\,e^{\pi i\,\s\{ \CM
(x\,\Om\,{}^tx+\,2\,x\,{}^tZ)\} },\quad (\Om,Z)\in\BH_{n,m},\
x\in\rmn.
\end{equation}

For brevity we put $\mfoz:=\mfm (\Om,Z)$ for $(\Om,Z)\in
\BH_{n,m}.$ We put
$$G_*^J:= G_*\ltimes H_\BR^{(n,m)}.$$
We observe that $G^J_*$ acts on $\BH_{n,m}$ through the natural projection of
$G^J_*$ onto $G^J$.
Let $J^*_\CM:G^J_*\times \BH_{n,m}\lrt \BC^{\times}$ be an automorphic
factor for $G^J_*$ on $\BH_{n,m}$ defined by

\begin{eqnarray}
   J_\CM^*(\widetilde g,(\Om,Z))&= e^{\pi
   i\,\sigma\left({\mathcal{M}}(Z+\lambda\,
\Om+\mu)(C\Om+D)^{-1}C\,{}^t(Z+\l\,\Om+\mu)\right)} \hskip 3cm \\
& \times\, e^{-\pi i\,\sigma \left( \mathcal{M}(\lambda\,
\Om\,^t\!\lambda\,+\,2\,\lambda\,{}^t\!Z\,+\,\kappa\,+\,
\mu\,{}^t\!\lambda ) \right) } J_{1/2}((g,\epsilon),\Om)^m,\nonumber
\end{eqnarray}
where ${\widetilde g}_*=((g,\epsilon),(\lambda,\mu;\kappa))\in G^J_*$ with $g=\begin{pmatrix} A&B\\
C&D\end{pmatrix}\in Sp(n,\BR),\ (\lambda,\mu;\kappa)\in
H_{\BR}^{(n,m)}$ and $(\Om,Z)\in \BH_{n,m}.$

\begin{theorem} Let $m$ be an odd positive integer. The map ${\mathscr F}^{(\CM)}:\BH_{n,m}\lrt
L^2\big(\rmn\big)$ defined by (7.1) is a covariant map for the
Schr{\"o}dinger-Weil representation $\omega_\CM$ of $G^J$ and the
automorphic factor $J^*_\CM$ for $G^J_*$ on $\BH_{n,m}$ defined by
Formula (7.2). In other words, $\mfm$ satisfies the following
covariance relation

\begin{equation}
\omega_\CM ({\widetilde g}_*)\mfoz=J_\CM^* \big( {\widetilde
g}_*,(\Om,Z)\big)^{-1} \mfm_{{\widetilde g}_*\cdot (\Om,Z)}
\end{equation}

\noindent for all ${\widetilde g}_*=((g,\epsilon),(\l,\mu;\k)\in G^J_*$
and $(\Om,Z)\in
\BH_{n,m}.$
\end{theorem}

\noindent {\it Proof.} The proof can be found in \cite{Y18}\ (cf.\,\cite{Y17}).
\hfill $\square$

\vskip 0.2cm For a positive definite integral matrix $\CM$ of degree $m$, we define the holomorphic function
$\Theta_\CM:\BH_{n,m}\lrt\BC$ by

\begin{equation}
\Theta_\CM (\Om,Z)=\sum_{\xi\in \BZ^{(m,n)}}
e^{\pi\,i\,\sigma\left( \CM(
\xi\,\Om\,{}^t\xi\,+\,2\,\xi\,{}^tZ)\right) },\quad (\Om,Z)\in
\BH_{n,m}.
\end{equation}

We can prove the following theorem.

\begin{theorem} The function $\Theta_\CM$ is a Jacobi form of weight ${\frac m2}$
and index ${{\CM}\over 2}$
with respect to a discrete subgroup $\G_{\CM,*}^J:=\G_{\CM,*} \ltimes
H_\BZ^{(n,m)}$ of $\G^J_*$ for a suitable arithmetic subgroup
$\G_\CM$ of $\G_n$ with $\G_{\CM,*}=\pi_*^{-1}(\G_\CM).$
That is, $\Theta_\CM$ satisfies the functional
equation
\begin{equation}
\Theta_\CM\big( {\widetilde \g}_*\cdot (\Om,Z)\big)=\rho_\CM (\g_*)\,
J^*_\CM\big( {\widetilde \g}_*,(\Om,Z)\big) \Theta_\CM(\Om,Z),\quad (\Om,Z)\in
\BH_{n,m},
\end{equation}
\end{theorem}

\noindent where $\rho_\CM$ is a suitable character of
$\G_{\CM,*}$ and ${\widetilde \g}_*=(\g_*,(\l,\mu;\k))\in \G_{\CM,*}^J$.

\noindent {\it Proof.} The proof can be found in \cite{Y17} when
$\CM$ is unimodular and even integral. In the case $\CM$ is a symmetric positive integral
matrix of odd degree $m$ such that $\det\,(\CM)=\,1$ with a special arithmetic subgroup $\G_{\CM,*}$,
the proof can be found in \cite{Y18}.
In a similar way we can
prove the above theorem. \hfill $\square$

\vskip 0.2cm According to Theorem 1 and Theorem 2, we see that the
theta series $\Theta_\CM$ is closely related to the
Schr{\"o}dinger-Weil representation of the Jacobi group $G^J$. We
note that the theta series
\begin{equation}
\Theta (\Om)= \,\sum_{A\in \BZ^n} e^{\pi\,i\,\sigma
(A\,\Om\,{}^t\!A)},\quad \ \Om\in \BH_n
\end{equation}

\noindent is a Siegel modular form of weight ${\frac 12}$ with
respect to the theta subgroup $\G_{\Theta}$ of $\G_n$, that is,
$\Theta$ satisfies the following functional equation
\begin{equation}
\Theta (\g\cdot\Om)=\,\zeta (\g)\,\left(\det
(C\Om+D)\right)^{\frac 12}\,\Theta(\Om),\quad \ \Om\in \BH_n,
\end{equation}

\noindent where $\zeta (\g)$ is a character of $\G_\Theta$ with
$|\zeta(\g)|^8=1$ and $\g=\begin{pmatrix} A & B\\ C & D
\end{pmatrix}\in \G_\Theta.$ We refer to \cite[pp.\,189-201]{Mum}
for more detail. Indeed the function ${\mathscr F}:\BH_n\lrt
L^2\big( \BR^n\big)$ defined by
\begin{equation}
{\mathscr F}(\Om)(x)=\,e^{\pi\,i\,\sigma (x\,\Om\,{}^tx)},\quad \
\Om\in \BH_n\ \textrm{and}\ x\in \BR^n.
\end{equation}

\noindent is a covariant map for the Weil representation $\omega$
of $Sp(n,\BR)$ and the automorphic form ${\mathfrak J}_{\frac
12}:Sp(n,\BR)\times \BH_n\lrt \BC^{\times}$ defined by
\begin{equation}
{\mathfrak J}_{\frac 12}(g,\Om)=\,\left(\det (C\Om+D)\right)^{\frac 12},\quad
\Om\in \BH_n
\end{equation}

\noindent with $g=\begin{pmatrix} A & B\\ C & D
\end{pmatrix}\in Sp(n,\BR).$ More precisely, if we put ${\mathscr
F}_\Om:=\,{\mathscr F}(\Om)$ for brevity, the vector valued map
${\mathscr F}$ satisfies the following covariance relation
\begin{equation}
\omega (g){\mathscr F}_{\Om}=\,\left(\det (C\Om+D)\right)^{-\frac
12}\,{\mathscr F}_{g\cdot\Om}
\end{equation}

\noindent for all $g\in Sp(n,\BR)$ and $\Om\in\BH_n.$ We refer to
\cite{LV} for more detail. This is a special case of Theorem 1 and
Theorem 2.

%%%%%%%%%%%%%%%%%%%%%%%%%%%%%%%%%%%%%%%%%%%%%%%%%%%%%%%%%%%%%%%%%%%%%%%%%%%%%%%%%%%%%%%%%%%%%%%%%%%%%%%%%%%%%%%%%%%%%%%%%%%%%%%%%%%%
%%%%%%%%%%%%%%%%%%%%%%%%%%%%%%%%%%%%%%%%%%%%%%%%%%%%%%%%%%%%%%%%%%%%%%%%%%%%%%%%%%%%%%%%%%%%%%%%%%%%%%%%%%%%%%%%%%%%%%%%%%%%%%%%%%%%

\vskip 0.885cm \noindent {\bf 7.2. Maass-Jacobi Forms }
\vskip 0.3cm
Recently in the case $n=m=1$ A. Pitale \cite{P} gave a new definition of nonholomorphic Maass-Jacobi forms of
weight $k$ and $m\in\BZ^+$ as eigenfunctions of a certain differential operator ${\mathcal C}^{k,m}$, and
constructed new examples of cuspidal Maass-Jacobi forms $F_f$ of even weight $k$ and index $1$
from Maass forms $f$ of weight half integral weight $k-1/2$ with respect to $\G_0(4)$. Moreover he also showed that
the map $f\mapsto F_f$ is Hecke equivariant and compatible with the representation theory of the Jacobi group $G^J.$
We will describe his results in some detail.

\vskip 0.153cm
For a positive integer $N$, we let
\begin{equation*}
\Gamma_0(N)=\,\left\{ \begin{pmatrix} a & b \\ c & d \end{pmatrix}\in SL(2,\BZ)\, \Big|\ c\equiv 0 \
( \textrm{mod}\ N)\ \right\}
\end{equation*}
be the congruence subgroup of $SL(2,\BZ)$ called the {\it Hecke subgroup} of level $N$.
Let ${\mathfrak G}$ be the group which consists of all pairs $(\gamma,\phi(\tau))$, where
$\gamma=\begin{pmatrix} a & b \\ c & d \end{pmatrix}\in GL(2,\BR)^+$ and $\phi(\tau)$ is a
function on ${\mathbb H}$ such that
\begin{equation*}
\phi(\tau)=\,t\,\det (\gamma)\,\left( {{(c\tau+d)}\over {|c\tau+d|}}\right)^{1/2}
\quad \textrm{with}\ t\in \BC,\ |t|=1.
\end{equation*}
The group law is given by
\begin{equation}
(\gamma_1,\phi_1(\tau))\cdot (\gamma_2,\phi_2(\tau))=\,\big(\gamma_1\gamma_2, \phi_1(\gamma_2\cdot\tau)\phi_2(\tau)\big),\quad
\gamma=\,\begin{pmatrix} a & b \\ c & d \end{pmatrix}\in GL(2,\BR)^+.
\end{equation}
Then there is an injective homomorphism $\G_0(4)\mapsto {\mathfrak G}$ given by
\begin{equation}
\g\mapsto \g^*:=\big(\g,j(\g,\tau) \big),
\end{equation}
where $\gamma=\,\begin{pmatrix} a & b \\ c & d \end{pmatrix}\in \G_0(4)$ and
$$j(\g,\tau):=\,\left( {\frac cd}\right)\,\epsilon_d^{-1}\,\left( {{(c\tau+d)}\over {|c\tau+d|}}\right)^{1/2}
=\,{{\theta(\g\cdot\tau)}\over {\theta(\tau)}}$$
with
$$ \theta(\tau):=\,y^{1/4}\,\sum_{n=-\infty}^{\infty} e^{2\pi i n^2 \tau}$$
and
\begin{equation*}
\epsilon_d = \begin{cases} 1,\quad \textrm{if}\ d\equiv 1\ ( \textrm{mod}\ 4),\\
i,\quad \textrm{if}\ d\equiv 3\ ( \textrm{mod}\ 4) \end{cases}.
\end{equation*}
And $\left( {\frac cd}\right)$ is defined as in \cite[p.\,442]{S}.

\vskip 0.2cm
For an integer $k\in\BZ$, we define the slash operator $||_{k-1/2}$ on functions on $\BH$ as follows\,:
\begin{equation}
\left( f||_{k-1/2}(\gamma,\phi)\right)(\tau):=\,f(\g\cdot \tau)\,\phi(\tau)^{-(2k-1)}.
\end{equation}

\vskip 0.2cm\noindent
\begin{definition}
A smooth function $f:\BH\lrt \BC$ is called a {\it Maass form} of weight $k-1/2$ with respect to $\G_0(4)$
if it satisfies the following properties (M1)-(M3)\,:
\vskip 0.1cm\noindent
\textrm{(M1)}\ \ $ f||_{k-1/2}\g^*=\,f\quad \textrm{for\ all}\ \g\in\G_0(4).$
\vskip 0.1cm\noindent
\textrm{(M2)}\ \ $\Delta_{k-1/2}f=\,\Lambda f \quad \textrm{for\ some}\ \Lambda\in\BC,\ \textrm{where}\
\Delta_{k-1/2}$ is the Laplace-Beltrami operator given by
\begin{equation}
\Delta_{k-1/2}= \,y^2\,\left( { {\partial^2\ \, }\over {\partial x^2}}+
{ {\partial^2\ \, }\over {\partial y^2}}
\right)-i \left( k-1/2\right) y { {\partial\ \, }\over {\partial x}}.
\end{equation}
\vskip 0.1cm\noindent
\textrm{(M3)}\ \ $f(\tau)=O(y^N)$\ as\ $y\lrt \infty$ \textrm{for some}\ $N>0.$
\vskip 0.1cm\noindent
If, in addition, $f$ vanishes at all the cusps of $\G_0(4)$, then we say that $f$ is a Maass cusp form.
\end{definition}

We denote by $M_{k-1/2}(4)$ (resp. $S_{k-1/2}(4)$) be the vector space of all Maass forms
(resp. Maass cusp forms) of weight $k-1/2$ with respect to $\G_0$. As shown in \cite{KP}
or \cite{P05}, if
$f\in M_{k-1/2}(4)$, then $f$ has the following Fourier expansion
\begin{equation}
f(\tau)=\sum_{n\in\BZ} c(n) W_{ \textrm{sgn} {{k-1/2}\over 2},{{il}\over 2} } (2\pi |n|y)\,e^{2\pi i nx},
\end{equation}
where $\Lambda=-\left\{ 1/4+(l/2)^2 \right\}$ and $W_{\mu,\nu}(y)$ is the classical Whittaker function which
is normalized so that $W_{\mu,\nu}(y)\sim e^{-y/2}y^{\mu}$ as $y\lrt\infty$. If $f\in S_{k-1/2}(4)$, then
we have $c(0)=0$ in (7.15). We define the plus space by
\begin{equation}
M_{k-1/2}^+(4):=\,\left\{ f \in M_{k-1/2}(4)\,|\ c(n)=0\ \textrm{whenever}\ (-1)^{k-1}n\equiv 2,3\
\textrm{mod} 4)\,\right\}.
\end{equation}
We set
\begin{equation*}
S_{k-1/2}^+(4):=M_{k-1/2}^+(4)\cap S_{k-1/2}(4).
\end{equation*}

For a given integer $k\in\BZ$ and $m\in \BZ^+$, we let
\begin{equation}
j_{k,m}^{ \textrm{nh}}({\tilde g},(\tau,z)):=\,e^{2\pi i\,m\left\{ \kappa-c(z+\l \tau+\mu)^2(c\tau+d)^{-1}
+\l^2\tau+2\l z+\l \mu\right\}} \times \left( {{c\tau+d}\over {|c\tau+d|}}\right)^{-k}
\end{equation}
be the nonholomorphic automorphic factor for $G^J$ on $\BH\times \BC$, where ${\tilde g}=(g,(\l,\mu;
\kappa))$ with $g=\begin{pmatrix} a & b \\ c & d \end{pmatrix}\in SL(2,\BR),\ \l,\mu,\kappa\in\BR$
and $(\tau,z)\in \BH\times\BC.$ For ${\tilde g}\in G^J(\BR),\ (\tau,z)\in \BH\times\BC$ and a smooth function
$F:\BH\times\BC\lrt\BC,$ we set
\begin{equation}
(F|_{k,m}{\tilde g})(\tau,z):=\,j_{k,m}^{ \textrm{nh}}({\tilde g},(\tau,z))F({\tilde g}\cdot (\tau,z)).
\end{equation}

\vskip 0.2cm
Let $\G^J:=SL(2,\BZ)\ltimes H_\BZ^{(1,1)}$ be the discrete subgroup of $G^J(\BR):=SL(2,\BR)\ltimes H_\BR^{(1,1)}.$

\vskip 0.2cm\noindent
\begin{definition}
A smooth function $F:\BH\times \BC\lrt \BC$ is called a {\it Maass-Jacobi form} of weight $k$ and
index $m$ with respect to $\G^J$
if it satisfies the following properties (MJ1)-(MJ3)\,:
\vskip 0.1cm\noindent
\textrm{(M1)}\ \ $ F({\tilde\g}\cdot (\tau,z))=\,j_{k,m}^{ \textrm{nh}}({\tilde \g},(\tau,z))^{-1}
F(\tau,z)\quad \textrm{for\ all}\ {\tilde\g}\in\G^J$ and $(\tau,z)\in \BH\times\BC$.
\vskip 0.1cm\noindent
\textrm{(M2)}\ \ ${\mathcal C}^{k,m}F=\,\lambda_{k,m}(f) F \quad \textrm{for\ some}\ \lambda_{k,m}(f)\in\BC.$
\vskip 0.1cm\noindent
\textrm{(M3)}\ \ $F(\tau,z)=O(y^N)$\ as\ $y\lrt \infty\ \textrm{for some}\ N>0$.
\vskip 0.1cm\noindent
If, in addition, $f$ satisfies the following cuspidal condition
\begin{equation}
\int_0^1\int_0^1 F\left( \begin{pmatrix} 1 & x \\ 0 & 1 \end{pmatrix} (0,\mu;0)(\tau,z)\right)\,
e^{-2\pi i(nx+r\mu)}dxd\mu=0
\end{equation}
for all $n,r\in\BZ$ such that $4mn-r^2=0$, then we say that $f$ is a Maass-Jacobi cusp form.
\vskip 0.2cm
In (M2), ${\mathcal C}^{k,m}$ is the $G^J(\BR)$-invariant differential operator defined by
\begin{eqnarray*}
{\mathcal C}^{k,m}F &=&\, {\frac 58}\,F-2(\tau-{\bar\tau})^2F_{\tau{\bar\tau}}-(k-1)(\tau-{\bar\tau})F_{\bar\tau}
-k(\tau-{\bar\tau})F_\tau \\
& &\ + { {k (\tau-{\bar\tau})} \over {8\pi i\,m} } F_{zz}\,+\,{ {(\tau-{\bar\tau})^2}\over {4\pi i\, m} }
F_{{\bar\tau}zz}
\, + \,{ {k\,(\tau-{\bar\tau})}\over {4\pi i\, m} }\,F_{z{\bar z}}\\
& &\ +\,
{ {(\tau-{\bar\tau})(z-{\bar z})}\over {4\,\pi i \,m}}\,F_{z z {\bar z}}\,-\,
2\,(\tau-{\bar\tau})(z-{\bar z})\,F_{\tau {\bar z}}\,+\, { {(\tau-{\bar\tau})^2}\over {4\,\pi i\,m}}\,
F_{\tau{\bar z}{\bar z}}\\
& & \ +\,\left( { {(z-{\bar z})^2}\over 2 }\,+\,{ {k(\tau-{\bar\tau})}\over {8\,\pi i\,m} }\right)
F_{{\bar z}{\bar z}}\,+\, { {(\tau-{\bar\tau})(z-{\bar z})}\over {4\,\pi i \,m}}\,
F_{z {\bar z} {\bar z}}.
\end{eqnarray*}
\end{definition}

\vskip 0.2cm\noindent
We denote by $J_{k,m}^{ \textrm{nh}}$ \big( resp. $J_{k,m}^{ \textrm{nh},cusp}$ \big) the vector space of all
Maass-Jacobi forms (resp. Maass-Jacobi cusp forms) of weight $k$ and index $m$ with respect to
$\G^J$.

\vskip 0.2cm For a Maass form $f\in M^+_{k-1/2}(4)$ with $k\in 2\BZ$, he defined the function $F_f$ on
$\BH\times \BC$ by
\begin{equation}
F_f(\tau,z):=\,f^{(0)}(\tau)\,{\widetilde{\Theta}}^{(0)}(\tau,z)\,+\,
f^{(1)}(\tau)\,{\widetilde{\Theta}}^{(1)}(\tau,z),\quad (\tau,z)\in\BH\times \BC.
\end{equation}
We refer to \cite[pp.\,96-97]{P} for the precise definition of $f^{(0)},\,f^{(1)},\,{\widetilde{\Theta}}^{(0)}$
and ${\widetilde{\Theta}}^{(1)}$. Pitale \cite{P} showed that if $f\in M^+_{k-1/2}(4)$ with $k\in 2\BZ$, then
$F_f\in J_{k,1}^{ \textrm{nh}}$, and $F_f\in J_{k,1}^{ \textrm{nh},cusp}$ if and only if $f\in S_{k-1/2}^+(4).$
Furthermore he showed that if $\Delta_{k-1/2}f=\Lambda f$, then ${\mathcal C}^{k,1}F_f=\,2\,\Lambda\, F_f$
under the assumption $f\in M^+_{k-1/2}(4)$ with $k\in 2\BZ$.

\vskip 0.2cm For an odd prime $p$, the Jacobi Hecke operator $T_p$ on $J_{k,1}^{ \textrm{nh}}$ (cf.
\cite[p.\,168]{BS} or \cite[p.\,41]{EZ} is defined by
\begin{equation}
T_pF:=\sum_{\substack{M\in SL(2,\BZ)/\BZ^{(2,2)}\\ \det(M)=p^2\\ \textrm{gcd}(M)=1} }
\sum_{(\l,\mu)\in (\BZ/p\BZ)^2}
F|_{k,1}\left( \det(M)^{-1/2}M(\l,\mu;0) \right).
\end{equation}

\vskip 0.2cm\noindent
\begin{theorem} Let $f\in S^+_{k-1/2}(4)\,(k\in 2\BZ)$ be a Hecke eigenform with eigenvalue $\l_p$ for every odd
prime $p$. Then $T_p=\,p^{k-3/2}\l_p\,F_f$ for all odd prime $p$. Namely $F_f$ is also an eigenfunction of all
$T_p$ for every odd prime $p$.
\end{theorem}
\noindent {\it Proof.} The proof can be found in \cite[pp.\,104-106]{P}. \hfill $\square$

\vskip 0.2cm Let $f$ be a Hecke eigenform in $S^+_{k-1/2}(4)\,(k\in 2\BZ)$ such that for every odd prime $p$ we have
$T_p f=\l_p f$ and $\Delta_{k-1/2}f=\Lambda f$ with $\Lambda={\frac 14}(s^2-1).$ Let ${\widetilde \pi}_f=\otimes
{\widetilde \pi}_{f,p}$ be the irreducible cuspidal genuine automorphic representation of a two-fold covering group
${\widetilde{SL(2,\BA)}}$ of $SL(2,\BA)$ corresponding to $f$ (cf.\,\cite[p.\,386]{W}). Now we let $F_f$ be the
Maass-Jacobi cusp form in $J_{k,1}^{ \textrm{nh},cusp}$ constructed from an eigenform $f\in S^+_{k-1/2}(4)\,(k\in 2\BZ)$
by Formula (7.20). Then $F_f$ is an eigenform of all $T_p$ for every odd prime $p$ and is an eigenfunction of
the differential operator ${\mathcal C}^{k,1}.$ We lift $F_f$ to the function $\Phi_{F_f}$ on $G^J(\BA)$
as follows. By the strong approximation theorem for $G^J(\BA)$, we have the decomposition
\begin{equation}
G^J(\BA)=\,G^J(\BZ)\, G^J(\BR)\, \Pi_{p< \infty} G^J(\BZ_p).
\end{equation}
If ${\tilde g}=\g {\tilde g}_\infty k_0\in G^J(\BA)$ with $\g\in G^J(\BZ),\, {\tilde g}_\infty\in G^J(\BR),\,
k_0\in \Pi_{p< \infty} G^J(\BZ_p),$ we define
\begin{equation}
\Phi_{F_f}({\tilde g}):=\,(F_f|_{k,m}{\tilde g}_\infty)(i,0).
\end{equation}
Let $\Pi_{F_f}$ be the space of all right translates of $\Phi_{F_f}$ on which $G^J(\BA)$ acts by right translation.
Pitale \cite{P} proved that
\begin{equation}
\Pi_{F_f}=\,{\widetilde \pi}_f\otimes \omega_{ \textrm{SW}}^1,
\end{equation}
where $\omega_{ \textrm{SW}}^1$ is the Schr{\"o}dinger-Weil representation of $G^J(\BA)$ (cf.\,\cite{BS}).

\begin{remark} For a Siegel cusp form of half integral weight, we have a result similar to Formula (7.24).
See \cite{BS} for the case $n=1$ and \cite{Ta3, Ta4} for the case $n\geq 1.$
\end{remark}

\begin{remark} Berndt and Schmidt \cite{BS} gave a definition of Maass-Jacobi forms different from Definition 7.2.
Yang \cite{Y12, Y14} gave a definition of Maass-Jacobi forms using the Laplacian of an invariant metric on the Siegel-Jacobi space
$\BH_n\times \BC^{(m,n)}$ in the aspect of the spectral theory on $L^2\big(\G_n^J\ba \BH_n\times\BC^{(m,n)}\big)$.
\end{remark}

%%%%%%%%%%%%%%%%%%%%%%%%%%%%%%%%%%%%%%%%%%%%%%%%%%%%%%%%%%%%%%%%%%%%%%%%%%%%%%%%%%%%%%%%%%%%%%%%%%%%%%%%%%%%%%%%%%%%%%%%%%%%%%%%%%%%%%
%%%%%%%%%%%%%%%%%%%%%%%%%%%%%%%%%%%%%%%%%%%%%%%%%%%%%%%%%%%%%%%%%%%%%%%%%%%%%%%%%%%%%%%%%%%%%%%%%%%%%%%%%%%%%%%%%%%%%%%%%%%%%%%%%%%%%%

\vskip 0.885cm \noindent {\bf 7.3. Theta Sums} \vskip 0.3cm

We embed $SL(2,\BR)$ into $Sp(n,\BR)$ by
\begin{equation}
SL(2,\BR)\ni \begin{pmatrix} a & b \\ c & d \end{pmatrix} \longmapsto
\begin{pmatrix} aI_n & bI_n \\ cI_n & d I_n\end{pmatrix}\in Sp(n,\BR).
\end{equation}
Every map $M=\begin{pmatrix} a & b \\ c & d \end{pmatrix}\in SL(2,\BR)$ admits the unique Iwasawa decomposition
\begin{equation*}
M=\begin{pmatrix} 1 & x \\ 0 & 1 \end{pmatrix}\begin{pmatrix} y^{1/2} & 0 \\ 0 & y^{-1/2} \end{pmatrix}
\begin{pmatrix} \cos \theta & -\sin\theta \\ \sin\theta & \ \ \cos\theta  \end{pmatrix}=\,(\tau,\theta),
\end{equation*}
where $\tau=x+\,i\,y\in\BH_1$ and $0\leq \theta < 2\pi.$ Then $SL(2,\BR)$ acts on $\BH_1\times [0,2\pi)$ by
\begin{equation}
M\cdot (\tau,\theta):=\,\big(M\cdot\tau, \,\theta+\arg(c\tau+d)\ \textrm{mod} \ 2\pi\big),
\end{equation}
where $M=\begin{pmatrix} a & b \\ c & d \end{pmatrix}\in SL(2,\BR),\ \tau\in\BH_1$ and $\theta\in [0,2\pi).$

\vskip 0.2cm
We put
\begin{equation*}
G_{n,1}^J:=\, Sp(n,\BR)\ltimes H_\BR^{(n,1)}.
\end{equation*}
We take ${\mathcal M}=1$ in Section 4. Then we let ${\mathscr W}=\,{\mathscr W}_{\mathcal M},
\ R=R_{\mathcal M}$ and $c=c_{\mathcal M}$\,(see Section 4).
If $M_i=\begin{pmatrix} a_i & b_i \\ c_i & d_i \end{pmatrix}\in SL(2,\BR)$ for $i=1,2,3$ with $M_3=M_1M_2$, then
the cocycle $c$ is given by
\begin{equation*}
c(M_1,M_2)=\,e^{-i\pi n \,\textrm{sign}(c_1c_2c_3)/4},
\end{equation*}
where
\begin{equation*}
\textrm{sign}(x)=\,\begin{cases} -1 \qquad \textrm{if}\ x<0,\\
\ \ 0 \qquad \textrm{if}\ x=0,\\
\ \ 1 \qquad \textrm{if}\ x<0.\end{cases}
\end{equation*}

For $(\tau,\theta)\in SL(2,\BR),$ we define
\begin{equation}
{\widetilde R}(\tau,\theta):=\,e^{-i\pi n\,\sigma_\theta /4}\,R(\tau,\theta),
\end{equation}
where
\begin{equation*}
\sigma_\theta   =\,\begin{cases} \ \ \ 2\nu \ \ \ \qquad \textrm{if}\ \theta=\nu\pi,\ \nu\in\BZ,\\
\ 2\nu+1 \qquad \textrm{if}\ \nu\pi <\theta < (\nu+1)\pi,\ \nu\in\BZ.
\end{cases}
\end{equation*}
Then ${\widetilde R}$ is a unitary representation of the double covering group of $SL(2,\BR)$ \ (cf.\,
\cite{LV}). Obviously ${\widetilde R}(i,\theta){\widetilde R}(i,\theta')=\,{\widetilde R}(i,\theta+\theta').$

\vskip 0.2cm
We see that
\begin{equation}
\omega_{ \textrm{SW}}^1((\xi;t)(\tau,\theta))=\,{\mathscr W}((\xi;t))\,{\widetilde R}(\tau,\theta),
\end{equation}
where $\omega_{ \textrm{SW}}^1$ denotes the Schr{\"o}dinger-Weil representation of $G_{n,1}^J$ (see
Formula (4.10)). Here $(\xi;t)\in H_\BR^{(n,1)}$ and $(\tau,\theta)$ is considered as an element of $Sp(n,\BR)$ by the
embedding (7.25).

\vskip 0.2cm
We denote by ${\mathcal S}(\BR^n)$ the vector space of $C^\infty$-functions on $\BR^n$ that, as well as their
derivatives, decrease rapidly at $\infty.$ For any $f\in {\mathcal S}(\BR^n)$,
{\it Jacobi's theta sum} for $f$ is defined to be the function
\begin{equation}
\Theta_f (\tau,\theta;\xi,t):=\,\sum_{\alpha\in\BZ^n}
\left[ \omega_{ \textrm{SW}}^1 \big( (\xi;t)(\tau,\theta)\big)f\right](\alpha),
\end{equation}
where $(\tau,\theta)\in SL(2,\BR)\hookrightarrow Sp(n,\BR)$ and $(\xi;t)\in H_\BR^{(n,1)}$ with
$\xi=(\lambda,\mu),\ \lambda,\mu\in \BR^n$ and $t\in\BR.$
For $f,g\in {\mathcal S}(\BR^n)$, the product of theta sums of the form
$$ \Theta_f (\tau,\theta;\xi,t)\,{\overline{ \Theta_g (\tau,\theta;\xi,t)}}$$
is independent of the $t$-variable.

\vskip 0.2cm Let us therefore define the semi-direct product group
$$G[n]:=\,SL(2,\BR)\ltimes \BR^{2n}$$
with multiplication law
$$(M,\xi)(M',\xi')=\,(MM',\,\xi+M\xi'),\quad M,M'\in SL(2,\BR),\ \xi,\xi'\in \BR^{2n}.$$
The set
\begin{equation*}
\Gamma [n]=: \left\{ \left( \begin{pmatrix} a & b \\ c & d \end{pmatrix},\,
\begin{pmatrix} a\,b\,{\mathfrak s} \\ c\,d\,{\mathfrak s} \end{pmatrix}+\alpha\right)\ \Big|\
\begin{pmatrix} a & b\\ c & d  \end{pmatrix}\in SL(2,\BZ),\ \alpha\in \BZ^{2n}\, \right\}
\end{equation*}
with ${\mathfrak s}=\,{}^t({\frac 12},{\frac 12},\cdots,{\frac 12})\in \BR^n$ is a subgroup of $G[n]$.
We can show that $\Gamma[n]$ is generated by
\begin{equation*}
\left( \begin{pmatrix} 0 & -1 \\ 1 & \ \ 0 \end{pmatrix},\,0\right),\ \
\left( \begin{pmatrix} 1 & 1 \\ 0 & 1 \end{pmatrix},
\begin{pmatrix} {\mathfrak s} \\ 0 \end{pmatrix} \right),\ \
\left(
\begin{pmatrix} 1 & 0\\ 0 & 1  \end{pmatrix},\, \alpha\right),\ \alpha\in \BZ^{2n}.
\end{equation*}

\noindent We put, for brevity,
$$\Theta_f (\tau,\theta;\xi):=\,\Theta_f (\tau,\theta;\xi,0).$$

J. Marklof \cite{Ma} proved the following properties of Jacobi's theta sums.

\begin{theorem} Let $f$ and $g$ be two elements in ${\mathcal S}(\BR^n)$. Then
\vskip 0.21cm
\noindent
(1) $\Theta_f (\tau,\theta;\xi)\,{\overline{ \Theta_g (\tau,\theta;\xi)}}$ is invariant under the
action of the left action of $\Gamma[n].$
\vskip 0.21cm
\noindent
(2) For any real number $R>1$, we have
\begin{eqnarray*}
& & \ \Theta_f (\tau,\theta;\xi)\,{\overline{ \Theta_g (\tau,\theta;\xi)}}\\
&=&\,y^{n/2}\,\sum_{\alpha\in \BZ^n} f_\theta \big( (\alpha-\mu)\,y^{1/2}\big)\,
{\overline{ g_\theta \big( (\alpha-\mu)\,y^{1/2}\big) }}\,+\,O_R \big(y^{-R}\big),
\end{eqnarray*}
where $\tau=x+\,i\,y\in \BH_1,\ \xi=(\lambda,\mu)$ with $\l,\mu\in\BR^n$ and
$$ f_\theta=\,{\widetilde R}(i,\theta)f.$$
\end{theorem}

\noindent
{\it Proof.} The proof can be found in \cite[pp.\,432-433]{Ma}.
\hfill $\square$

\vskip 0.3cm The above properties of Jacobi's theta sums together with Ratner's classification of
measures invariant under unipotent flows (cf.\,\cite{R1, R2}) are used to prove the important fact that under explicit
diophantine conditions on $(\alpha,\beta)\in \BR^2$, the local two-point correlations of
the sequence given by the values
$(m-\alpha)^2+(n-\beta)^2$ with $(m,n)\in\BZ^2$, are those of a Poisson process (see
\cite{Ma} for more detail).

\end{section}

\vskip 1cm
%%%%%%%%%%%%%%%%%%%%%%%%%%%%%%%%%%%%%%%%%%%%%%%%%%%%%%%%%%%%%%%%%%%%%%%%%%%%%%%%%%%%%%%%%%%%%%%%%%%%%%%%%%%%%%%%%%%%%%%%%%%%%%%%%%%%
%%%%%%%%%%%%%%%%%%%%%%%%%%%%%%%%%%%%%%%%%%%%%%%%%%%%%%%%%%%%%%%%%%%%%%%%%%%%%%%%%%%%%%%%%%%%%%%%%%%%%%%%%%%%%%%%%%%%%%%%%%%%%%%%%%%%
%%%%%%%%%%%%%%%%%%%%%%%%%%%%%%%%%%%%%%%%%%%%%%%%%%%%%%%%%%%%%%%%%%%%%%%%%%%%%%%%%%%%%%%%%%%%%%%%%%%%%%%%%%%%%%%%%%%%%%%%%%%%%%%%%%%%
%%
%%
%%     Section 8  Applications of the Weil-Satake Representation
%%
%%
%%%%%%%%%%%%%%%%%%%%%%%%%%%%%%%%%%%%%%%%%%%%%%%%%%%%%%%%%%%%%%%%%%%%%%%%%%%%%%%%%%%%%%%%%%%%%%%%%%%%%%%%%%%%%%%%%%%%%%%%%%%%%%%%%%%%
%%%%%%%%%%%%%%%%%%%%%%%%%%%%%%%%%%%%%%%%%%%%%%%%%%%%%%%%%%%%%%%%%%%%%%%%%%%%%%%%%%%%%%%%%%%%%%%%%%%%%%%%%%%%%%%%%%%%%%%%%%%%%%%%%%%%
%%%%%%%%%%%%%%%%%%%%%%%%%%%%%%%%%%%%%%%%%%%%%%%%%%%%%%%%%%%%%%%%%%%%%%%%%%%%%%%%%%%%%%%%%%%%%%%%%%%%%%%%%%%%%%%%%%%%%%%%%%%%%%%%%%%%

\begin{section}{{\large\bf Applications of the Weil-Satake Representation}}
\setcounter{equation}{0}

\vskip 0.21cm In this section we provide some applications of the
Weil-Satake Representation ${\widehat\omega}_\CM:={\widehat\omega}_{\CM,iI_n}$ to the
theory of representations of the Jacobi group $G^J$. Throughout
this section, for brevity, we put $G:=Sp(n,\BR)$. We will keep
the notations and the conventions in Section 5. We recall the
notations $G_2=G_{2,iI_n}$ and $G^J_2=G_2\ltimes H_\BR^{(n,m)}$ in
Section 5. For a real Lie group ${\mathfrak G}$, we denote by
$\widehat{\mathfrak G}$ the unitary dual of ${\mathfrak G}.$ We
define the following projections
\begin{eqnarray*}
&& p_2:G_2\lrt G,\quad\ \ \ \,\,(g,t)\longmapsto g,\\
&& p^J: G^J\lrt G, \quad\ \ \ \,\,(g,h)\longmapsto g, \\
&& p^J_2:G_2^J \lrt G^J,\quad\ \ \left( (g,t),h\right)\longmapsto (g,h),\\
&& p_{2,J}:G_2^J \lrt G_2,\quad \ \big( (g,t),h\big)\longmapsto
(g,t).
\end{eqnarray*}

\noindent Let ${\mathscr Z}$ be the center of $G^J$. Obviously
${\mathscr Z}\cong S(m).$
\begin{proposition} Let $\chi_\CM$ be the character of ${\mathscr
Z}$ defined by $\chi_\CM (\k)=e^{2\,\pi\,i\,\sigma (\CM\k)}$ with
$\k\in {\mathscr Z}.$ We denote by $\widehat{G_2^J}\big(
{\overline\chi}_\CM\big)$ the set of all equivalence classes of
irreducible representations $\eta$ of $G_2^J$ such that $\eta
(\k)=\,\chi_\CM (\k)^{-1}$ for all $\k\in {\mathscr Z}.$ We put
${\widetilde \pi}=\pi\circ p_{2,J}$ for any $\pi\in
\widehat{G_2^J}.$ Then the correspondence
\begin{equation*}
\widehat{G_2}\lrt \widehat{G_2^J}\big(
{\overline\chi}_\CM\big),\quad\ \pi\longmapsto {\widetilde \pi}
\otimes {\widehat\omega}_\CM
\end{equation*}

\noindent is a bijection from $\widehat{G_2}$ to
$\widehat{G_2^J}\big( {\overline\chi}_\CM\big)$. Furthermore $\pi$
is square integrable if and only if ${\widetilde \pi} \otimes
{\widehat\omega}_\CM$ is square square integrable modulo
${\mathscr Z}.$
\end{proposition}

\noindent {\it Proof.} The proof can be found in \cite{Ta1}.
\hfill $\square$

\vskip 0.2cm We now consider a holomorphic discrete series
representation of $G^J$. Let $K$ be the stabilizer of the action
(1.1) at $iI_n$. Then
\begin{equation*}
K=\left\{ \,\begin{pmatrix} A & -B \\ B & \ A\end{pmatrix}\in G\
\Big|\ \ A\,+\,i\,B \in U(n)\ \right\}.
\end{equation*}

\noindent Thus $K$ can be identified with the unitary group
$U(n)$. Let $(\rho,V_\rho)$ be an irreducible representation of
$K$ with highest weight $\rho=(\rho_1,\cdots,\rho_n)\in\BZ^n$ such
that $\rho_1\geq \cdots\geq \rho_n\geq 0.$ Then $\rho$ can be
extended to a rational representation of $GL(n,\BC)$ that is also
denoted by $\rho$. The representation space $V_\rho$ of $\rho$ has
a Hermitian inner product $\langle\ \,,\ \rangle$ on $V_\rho$ such
that $\langle \rho(g)u,v\rangle=\langle v,\rho(g^*)v\rangle$ for
all $g\in GL(n,\BC),\ u,v\in V_\rho,$ where $g^*=\,{}^t{\overline
g}.$ We define the unitary representation $\tau_\rho$ of $K$ by
\begin{equation}
\tau_\rho (k):=\,\rho\big( J(k,iI_n)\big),\quad k\in K.
\end{equation}

\noindent For all ${\widetilde g}=(g,h)\in G^J$ with $g\in G$ and
$(\Om,Z)\in\BH_{n,m},$ we define
\begin{equation}
J_{\rho,\CM}\big(
\wg,(\Om,Z)\big):=\,J_\CM(\wg,(\Om,Z))\,\rho(J(g,\Om)).\quad
(\textrm{see}\ (5.1)\ \textrm{and}\ (5.2))
\end{equation}

\noindent We note that for all $\wg\in G^J,\ (\Om,Z)\in \BH_{n,m}$
and $u,v\in V_\rho$, we have the relation
\begin{equation*}
\langle \,J_{\rho,\CM}\big(
\wg,(\Om,Z)\big)u,v\,\rangle=\,\langle\, u, J_{\rho,\CM}\big(
\wg,(\Om,Z)\big)^*v \,\rangle,
\end{equation*}

\noindent where
\begin{equation*}
J_{\rho,\CM}\big( \wg,(\Om,Z)\big)^*=\, \overline{J_{\CM}\big(
\wg,(\Om,Z)\big)}\,\rho\big( \,{}^t{\overline{J(g,\Om)}}\big).
\end{equation*}

\noindent We let ${\mathbb E}_{\rho,\CM}$ be the Hilbert space
consisting of $V_\rho$-valued measurable functions $f$ on
$\BH_{n,m}$ satisfying the condtion
\begin{equation*}
(f,f)=||f||^2=\int_{\BH_{n,m}}\langle\,\rho(Y)f(\Om,Z),f(\Om,Z)\,\rangle\,\k_\CM(\Om,Z)\,dv,
\end{equation*}

\noindent where $\k_\CM (\Om,Z)$ and $dv$ are defined in (5.3) and
(6.15) respectively. We let $K^J:=K\times S(m)$ be a
subgroup of $G^J.$ The representation $\Pi_{\rho,\CM}:=
\textrm{Ind}_{K^J}^{G^J}(\rho\otimes \overline{\chi}_\CM )$
induced from a representation $\rho\otimes \overline{\chi}_\CM$ is
realized on ${\mathbb E}_{\rho,\CM}$ as follows: for any $\wg\in
G^J$ and $f\in {\mathbb E}_{\rho,\CM}$, $\Pi_{\rho,\CM}$ is given
by
\begin{equation}
\big( \Pi_{\rho,\CM}(\,\wg\,)f\big) (\Om,Z)=\,J_{\rho,\CM}(\,\wg,(\Om,Z))^{-1}f(\,
\wg^{-1}\cdot(\Om,Z)).
\end{equation}

\noindent Let ${\mathbb H}_{\rho,\CM}$ be the subspace of
${\mathbb E}_{\rho,\CM}$ consisting of holomorphic functions in
${\mathbb E}_{\rho,\CM}$. It is easily seen that ${\mathbb
H}_{\rho,\CM}$ is a closed subspace of ${\mathbb E}_{\rho,\CM}$
invariant under the action of $\Pi_{\rho,\CM}.$ We let
$\pi_{\rho,\CM}$ be the restriction of $\Pi_{\rho,\CM}$ to
${\mathbb H}_{\rho,\CM}$.

\vskip 0.15cm Takase \cite{Ta2} proved the following result.
\begin{theorem} Suppose $\rho_n > n+{\frac m2}.$ Then ${\mathbb H}_{\rho,\CM}\neq
0$ and $\pi_{\rho,\CM}$ is an irreducible representation of $G^J$
which is square integrable modulo ${\mathscr Z}.$ Moreover the
multiplicity of $\rho$ in the restriction $\pi_{\rho,\CM}|_K$ of
$\pi_{\rho,\CM}$ to $K$ is equal to one.
\end{theorem}

We let
\begin{equation*}
K_2=p_2^{-1}(K)=\,\left\{\,(k,t)\in K\times T\,|\ \,t^2=\,\det
J(k,iI_n)\,\right\}.
\end{equation*}

\noindent The Lie algebra ${\mathfrak k}$ of $K_2$ and its Cartan
subalgebra ${\mathfrak h}$ are given by
$${\frak k}=\left\{\begin{pmatrix} A&-B\\  B& \ A
\end{pmatrix} \in \mathbb R^{(2n,2n)}\,\bigg|\
A+\,^tA=0,\ B=\,^tB\right\}$$ and
$${\frak h}=\left\{\begin{pmatrix} 0 & -C\\ C & \ \,0
\end{pmatrix} \in \mathbb R^{(2n,2n)}\,\bigg|\
C=\text{diag}\,(c_1,c_2,\cdots,c_n)\ \right\}.$$ Here
$\text{diag}\,(c_1,c_2,\cdots,c_n)$ denotes the diagonal matrix of
degree $n$. We define $\lambda_j\in {\frak h}_{\mathbb C}^*$ by
$\lambda_j\left(\begin{pmatrix} 0 & -C\\ C &
0\end{pmatrix}\right): =\sqrt{-1}\,c_j.$ We put
$${\mathbb M}^+=\left\{\,\sum_{j=1}^n m_j\lambda_j\,\,\bigg|\ \,m_j\in {\frac 12}\mathbb Z,\ m_1\geq
\cdots\geq m_n,\ m_i-m_j\in \mathbb Z\ \text{for\ all}\
i,j\right\}.$$ We take an element $\lambda=\sum_{j=1}^n
m_j\lambda_j\in {\mathbb M}^+.$ Let $\tau$ be an irreducible
representation of $K$ with highest weight
$\tau=(\tau_1,\cdots,\tau_n) \in \mathbb Z^n,$ where
$\tau_j=m_j-m_n\,(1\leq j\leq n-1).$ Let $\tau_{[\lambda]}$ be the
irreducible representation of $K_2$ defined by
\begin{equation}
\tau_{[\lambda]}(k,t):=\,t^{2m_n}\cdot \tau(J(k,iI_n)),\ \
(k,t)\in K_2.
\end{equation}
Then $\tau_{[\lambda]}$ is the irreducible representation of $K_2$
with highest weight $\lambda=(m_1,\cdots,m_n)$ and $\lambda
\longmapsto\tau_{[\lambda]}$ is a bijection from ${\mathbb M}^+$
to ${\widehat K}_2,$ the unitary dual of $K_2.$ According to
\cite[Theorem 7.2]{KV}, we have a decomposition of the restriction
$\widehat{\omega}_\CM|_{K_2}$ into irreducible components\,:
$$\widehat{\omega}_\CM |_{K_2}=\bigoplus_{\lambda}\,m_\lambda\,\tau_{[\lambda]},$$
where $\lambda$ runs over
\begin{align*}
\ \ \ &\lambda=\sum_{j=1}^{s}\tau_j\lambda_j+{\frac
m2}\sum_{j=1}^n\lambda_j\in {\mathbb M}^+\
\ (s=\textrm{Min}\,\{m,n\}),\\
\ \ \ &\tau_j\in \mathbb Z\ \,\textrm{such\ that}\
\tau_1\geq\tau_2\geq\cdots\geq \tau_{s}\geq 0
\end{align*}
and the multiplicity $m_\lambda$ is given by
$$m_\lambda=\prod_{1\leq i< j \leq m}\left( 1+{{\tau_i-\tau_j}\over {j-i}}\right),$$
where $\tau_j=0$ if $j>s.$ Let ${\widehat G}_{2,d}$ be the set of
all equivalence classes of square integrable irreducible unitary
representations of $G_2.$ The correspondence
$$\pi\longmapsto \textrm{Harish-Chandra\ parameter\ of}\ \pi$$
is a bijection from ${\widehat G}_{2,d}$ to $\Lambda^+,$ where
$$\Lambda^+=\left\{\sum_{j=1}^n m_j\lambda_j\in {\mathbb M}^+\ \bigg|\ m_1>\cdots>m_n,\
m_i-m_j\neq 0\ \text{for\ all}\ i,j,\ i\neq j\right\}.$$ See
\cite{Wa}, Theorem 10.2.4.1 for the details.
\smallskip

We choose an element $\lambda=\sum_{j=1}^n m_j\,\lambda_j\in
{\mathbb M}^+.$ Let $\pi^{\lambda}\in {\widehat G}_{2,d}$ be the
representation corresponding to the Harish-Chandra parameter
$$\sum_{j=1}^n(m_j-j)\,\lambda_j\in \Lambda^+.$$
The representation $\pi^{\lambda}$ is realized  as follows\,(see
\cite{Kn}, Theorem 6.6)\,: Let $(\tau,V_{\tau})$ be the
irreducible representation of $K$ with highest weight
$\tau=(\tau_1,\cdots,\tau_n),\ \tau_i=m_i-m_n\,(\,1\leq j\leq
n-1\,).$ Let ${\mathscr H}^{\lambda}$ be a Hilbert space
consisting of $V_{\tau}$-valued holomorphic functions $\varphi$ on
$\BH_n$ such that
$$|\varphi|^2=\int_{\BH_n}\left(\tau (Y)\,
\varphi(\Om),\varphi(\Om)\right)\,(\det Y)^{m_n}\,dv_\Om <
\infty,$$ where $dv_\Om=\,(\det Y)^{-(n+1)}[dX]\wedge [dY]$ is a
$G$-invariant volume element on $\BH_n$. Then $\pi^{\lambda}$ is
realized on ${\mathscr H}^{\lambda}$ as follows\,: for any
$\sigma=(g,t)\in G_2$ and $f\in {\mathscr H}^{\lambda}$,
\begin{equation}
\left(\pi^{\lambda}(\sigma)f\right)(\Om)=J_{[\lambda]}(\sigma^{-1},\Om)^{-1}
f (\sigma^{-1} \Omega)
\end{equation}

\noindent for all $\sigma=(g,t)\in G_2$ and $ f\in {\mathscr
H}^{\lambda}.$ Here
\begin{equation*}
J_{[\lambda]}(\sigma,\Om)=\,\left\{\,
t\,\beta_{iI_n}(g,g^{-1})\,|\det J(g,\Om)|^{\frac 12}\, { {\g
(g\Om,g(iI_n))} \over {\g (\Om,iI_n)}
}\right\}^{m_n}\,\tau(J(g,\Om)).
\end{equation*}

\begin{theorem}
 Suppose $\tau_n>n+{\frac
m2}.$ We put $\lambda=\sum_{j=1}^n(\tau_j-{\frac m2})\lambda_j\in
{\mathbb M}^+.$ Then the unitary representation
$\pi_{\tau,\CM}\circ p_2^J$ of $G_2^J$ is unitarily equivalent to
the representation $(\pi^{\lambda}\circ p_{2,J})\otimes
{\widehat\omega}_\CM.$
\end{theorem}

\noindent {\it Proof.} The proof can be found in \cite{Ta2}.
\hfill $\square$

\vskip 0.321cm
Using Theorem 8.2, Takase \cite{Ta4} established a bijective correspondence between the space of cuspidal
Jacobi forms and the space of Siegel cusp forms of half integral weight which is compatible with the action of
Hecke operators. For example, the classical result (cf.\, \cite{EZ} and \cite{IB})
\begin{equation}
J_{k,1}^{   \textrm{cusp}}(\G_n)\cong S_{k-1/2}(\G_0(4))
\end{equation}
can be obtained by the method of the representation theory. Here $\G_n$ denotes the Siegel modular group of degree $n$
and $\G_0(4)$ denotes the Hecke subgroup of $\G_n$.

\end{section}

%%%%%%%%%%%%%%%%%%%%%%%%%%%%%%%%%%%%%%%%%%%%%%%%%%%%%%%%%%%%%%%%%%%%%%%%%%%%%%%%%%%%%%%%%%%%%%%%%%%%%%%%%%%%%%%%%%%%%%%%%%%%%%%%%%%%
%%%%%%%%%%%%%%%%%%%%%%%%%%%%%%%%%%%%%%%%%%%%%%%%%%%%%%%%%%%%%%%%%%%%%%%%%%%%%%%%%%%%%%%%%%%%%%%%%%%%%%%%%%%%%%%%%%%%%%%%%%%%%%%%%%%%
%%%%%%%%%%%%%%%%%%%%%%%%%%%%%%%%%%%%%%%%%%%%%%%%%%%%%%%%%%%%%%%%%%%%%%%%%%%%%%%%%%%%%%%%%%%%%%%%%%%%%%%%%%%%%%%%%%%%%%%%%%%%%%%%%%%%
%%
%%
%%           Section 9  Lifting of Automorphic Forms
%%
%%%%%%%%%%%%%%%%%%%%%%%%%%%%%%%%%%%%%%%%%%%%%%%%%%%%%%%%%%%%%%%%%%%%%%%%%%%%%%%%%%%%%%%%%%%%%%%%%%%%%%%%%%%%%%%%%%%%%%%%%%%%%%%%%%%%
%%%%%%%%%%%%%%%%%%%%%%%%%%%%%%%%%%%%%%%%%%%%%%%%%%%%%%%%%%%%%%%%%%%%%%%%%%%%%%%%%%%%%%%%%%%%%%%%%%%%%%%%%%%%%%%%%%%%%%%%%%%%%%%%%%%%
%%%%%%%%%%%%%%%%%%%%%%%%%%%%%%%%%%%%%%%%%%%%%%%%%%%%%%%%%%%%%%%%%%%%%%%%%%%%%%%%%%%%%%%%%%%%%%%%%%%%%%%%%%%%%%%%%%%%%%%%%%%%%%%%%%%%

%\begin{section}{{\bf Connections to Lifting of Automorphic Forms}}
%\setcounter{equation}{0}

%\vskip 0.21cm

%\end{section}

\vskip 1cm
\bibliography{central}

\begin{thebibliography}{3}
\parskip=0pt
\itemsep=0pt

\bibitem{AAG} S. T. Ali, J. P. Antoine and J.-P. Gazeau,
{\em Coherent states, wavelets, and their generalizations},
Springer-Verlag, New York (2000).




\bibitem{B1} S. Berceau, {\em Coherent states associated to the Jacobi group-a variation on
a theme by Erich K{\"a}hler},
Jour. of Geometry and Symmetry in Physics, Vol. {\bf 9} (2007), 1-8.

\bibitem{B2} S. Berceau, {\em Coherent states associated to the Jacobi group},
Romanian Reports in Physics, Vol. {\bf 59}, No. 4 (2007), 1089--1101.

\bibitem{B3} S. Berceau, {\em Generalzed squeezed states for the Jacobi group},
arXiv:0812.0717v1 [math.DG] 3 Dec 2008.


\bibitem{BS} R. Berndt and R. Schmidt, {\em Elements of the
Representation Theory of the Jacobi Group}, Birkh{\"a}user, 1998.

\bibitem{BCR} K. Bringmann, C. Conley and O. Richter,
{\em Maass Jacobi forms over complex quadratic fields}, Math. Res. Lett. $\textbf{14}$, No. 1 (2007),
137--156.

\bibitem{EZ} M. Eichler and D. Zagier, {\em The Theory of Jacobi Forms}, Progress in Math., $ \textbf{55}$,
Birkh{\"a}user, Boston, Basel and Stuttgart, 1985.

\bibitem{F} E. Freitag, {\em Siegelsche Modulfunktionen}, Grundlehren de mathematischen Wissenschaften {\bf 55},
Springer-Verlag, Berlin-Heidelberg-New York (1983).


\bibitem{Ge} S. Gelbart, {\em Weil's Representation and the Spectrum of the Metaplectic
Group}, Lecture Notes in Math. $ \textbf{530}$, Springer-Verlag,
Berlin and New York, 1976.

\bibitem{GS1} V. Guillemin and S. Sternberg, {\em Geometric Asymptotics}, Amer. Math. Soc.,
Providence, R. I. (1977).

\bibitem{GS2} V. Guillemin and S. Sternberg, {\em Symplectic Technique in Physics},
Cambridge University Press, Cambridge, 1984.

\bibitem{H} J. N. Hollenhors, Phys. Rev. D, {\bf 19} (1979), 1669--1679.


\bibitem{IB} T. Ibukiyama, {\em On Jacobi forms and Siegel modular forms of half integral weights},
Comment. Math. Univ. Sancti Pauli {\bf 41} (1992), 109-124.


\bibitem{I} T. Ikeda, {\em On the lifting of elliptic cusp forms to Siegel cusp forms of degree $2n$},
Ann. Math. {\bf 154} (2001), 641--681.


\bibitem{KV}  M. Kashiwara and M. Vergne, {\em On the
Segal-Shale-Weil Representations and Harmonic Polynomials},
Invent. Math. $ \textbf{44}$ (1978), 1--47.

\bibitem{KP} S. Katok and P. Sarnak, {\em Heegner points, cycles and Maass forms},
Israel Math. J. {\bf 84} (1993), 237--268.




\bibitem{Kn} A. W. Knapp, {\em Representation Theory of
Semisimple Groups}, Princeton Univ. Press, Princeton, (1986).

\bibitem{KS1} S. Kudla and S. Rallis, {\em On the Weil-Siegel formula},
J. reine angew. Math. {\bf 387} (1988), 1--68.

\bibitem{KS2} S. Kudla and S. Rallis, {\em On the Weil-Siegel formula II},
J. reine angew. Math. {\bf 391} (1988), 65-84.

\bibitem{KS3} S. Kudla and S. Rallis, {\em A regularized Weil-Siegel formula : the first term identity},
Ann. of Math.(2) {\bf 140} (1994), no. 1, 1-80.

\bibitem{LV}  G. Lion and M. Vergne, {\em  The Weil representation, maslov index and Theta seires},
Progress in Math., $ \textbf{6}$, Birkh{\"a}user, Boston, Basel
and Stuttgart, 1980.


\bibitem{Lu} E. Y. C. Lu, Lett. Nuovo. Cimento {\bf 2} (1971), 1241--1244.



\bibitem{Ma} J. Marklof, {\em Pair correlation densities of inhomogeneous quadratic forms},
Ann. Math. {\bf 158} (2003), 419-471.


\bibitem{Mum} D. Mumford, \textit{Tata Lectures on Theta I,} Progress
in Math. {\bf 28}, Boston-Basel-Stuttgart (1983).

\bibitem{Ni}  S. Niwa, {\em Modular forms of half-integral weight and the integral of certain theta-functions},
Nagoya Math. J., {\bf 56}\,(1975), 147--161.

\bibitem{Per} A. M. Perelomov, {\em Generalized Coherent states and their Applications},
Springer, Berlin (1986).


\bibitem{P05} A. Pitale, {\em Lifting from ${\widetilde {SL_2}}$ to $GSpin(1,4)$},
Int. Math. Res. Not. $ \textbf{63}$ (2005), 3919--3966.


\bibitem{P} A. Pitale, {\em Jacobi Maass forms}, Abh. Math. Sem. Hamburg $ \textbf{79}$ (2009), 87--111.



%\bibitem{Sa} M. Saito, {\em Representations unitaires des groupes symplectiques},
%J. Math. Soc. Japan, {\bf 24}, No. 2\,(1972), 232--251.


\bibitem{R1} M. Ratner, {\em On Raghunathan's measure conjecture  }, Ann. of Math., {\bf 134}
(1991), 545--607.

\bibitem{R2} M. Ratner, {\em Raghunathan's topological conjecture and distributions of
unipotent flows,}, Duke Math. J. {\bf 63}
(1991), 235--280.

\bibitem{Sat} I. Satake, {\em Fock representations and theta functions},
Ann. Math. Study {\bf 66}\,(1969), 393--405.



\bibitem{S} G. Shimura, {\em On modular forms of half integral
weight }, Ann. of Math., {\bf 97}(1973), 440--481.


\bibitem{Sh}  T. Shintani, {\em On construction of holomorphic cusp forms of half integral weight},
Nagoya Math. J., {\bf 58}\,(1975), 83--126.

\bibitem{Si}  C. L. Siegel, {\em Indefinite quadratische Formen und Funnktionentheorie I and II},
Math. Ann. {\bf 124}\,(1951), 17--54 and Math. Ann. {\bf
124}\,(1952), 364--387\,; Gesammelte Abhandlungen, Band III,
Springer-Verlag (1966), 105--142 and 154--177.

\bibitem{St} P. Stoler, Phys. Rev. D. {\bf 1} (1970), 3217--3219.


\bibitem{Ta1} K. Takase, {\em A note on automorphic forms}, J. reine angew. Math., {\bf 409}(1990),
138-171.

\bibitem{Ta2} K. Takase, {\em On unitary representations of Jacobi groups},
J. reine angew. Math., {\bf 430}(1992), 130-149.

\bibitem{Ta3} K. Takase, {\em On Two-fold Covering Group of $Sp(n,\BR)$ and Automorphic Factor of Weight 1/2},
Comment. Math. Univ. Sancti Pauli {\bf 45} (1996), 117-145.

\bibitem{Ta4} K. Takase, {\em On Siegel Modular Forms of Half-integral Weights and Jacobi Forms},
Trans. of American Math. Soc. $ \textbf{351}$, No. 2 (1999), 735--780.

\bibitem{W} J.-L. Waldspurger, {\em Sur les coefficients de Fourier de formes modulaires de poids demi-entier},
J. Math. Pures Appl. {\bf 60} (1981), 375--484.

\bibitem{Wa} G. Warner, Harmonic Analysis on Semisimple Lie
Groups, I, II, Springer-Verlag, Berlin-Heidelberg-New York,
 (1972).


\bibitem{W} A. Weil, {\em Sur certains groupes d'operateurs unitares},
Acta Math., {\bf 111}\,(1964), 143--211.

\bibitem{Wo} K. B. Wolf, {\em Geometric Optics on Phase Space}, Springer (2004).


\bibitem{Y1} J.-H. Yang, {\em Harmonic Analysis on the Quotient Spaces of
Heisenberg Groups}, Nagoya Math. J., {\bf 123}(1991), 103--117.

\bibitem{Y2} J.-H. Yang, {\em Harmonic Analysis on the Quotient Spaces of
Heisenberg Groups II }, J. Number Theory, {\bf 49}(1)(1994),
63--72.

\bibitem{Y3} J.-H. Yang, {\em A decomposition theorem on
differential polynomials of theta functions of high level},
Japanese J. of Mathematics, The Mathematical Society of Japan, New
Series, {\bf 22}(1)(1996), 37--49.

\bibitem{Y4} J.-H. Yang, {\em Fock Representations of the Heisenberg Group
$H_{\BR}^{(g,h)}$}, J. Korean Math. Soc., {\bf 34}, no. 2\,(1997),
345--370.


\bibitem{Y5} J.-H. Yang, {\em Lattice Representations of the Heisenberg Group
$H_{\BR}^{(g,h)}$}, Math. Annalen, {\bf 317}(2000), 309--323.


\bibitem{Y6}  J.-H. Yang, {\em The Siegel-Jacobi Operator},
Abh. Math. Sem. Univ. Hamburg $ \textbf{63}$ (1993), 135--146.

\bibitem{Y7}  J.-H. Yang, {\em Remarks on Jacobi forms of higher
degree}, Proc. of the 1993 Workshop on Automorphic Forms and
Related Topics, the Pyungsan Institute for Mathematical Sciences,
Seoul (1993), 33--58.

\bibitem{Y8} J.-H. Yang, {\em Singular Jacobi Forms}, Trans.
of American Math. Soc. $ \textbf{347}$, No. 6 (1995), 2041--2049.

\bibitem{Y9} J.-H. Yang, {\em Construction of Modular Forms from
Jacobi Forms}, Canadian J. of Math. $ \textbf{47}$ (1995),
1329--1339.

\bibitem{Y10} J.-H. Yang, {\em A geometrical theory of Jacobi
forms of higher degree}, Proceedings of Symposium on Hodge Theory
and Algebraic Geometry\,(\,edited by Tadao Oda\,), Sendai, Japan
(1996), 125-147 {\it or} Kyungpook Math. J. $ \textbf{40\,(2)}$,
209--237 (2000) {\it or} arXiv:math.NT/0602267.

\bibitem{Y11} J.-H. Yang, {\em The method of orbits for real Lie groups},
Kyungpook Math. J. $ \textbf{42}$ (2002), 199--272.

\bibitem{Y12} J.-H. Yang, {\em A Note on Maass-Jacobi Forms},
Kyungpook Math. J. $ \textbf{43}$ (2003), 547--566.

\bibitem{Y13} J.-H. Yang, \textit{A note on a fundamental domain for Siegel-Jacobi space,} Houston Journal of
Mathematics, Vol. $ \textbf{32}$, No. 3 (2006), 701--712.

\bibitem{Y14} J.-H. Yang, {\em Invariant metrics and Laplacians on Siegel-Jacobi space,} arXiv:math.NT/0507215 v1 or
Journal of Number Theory {\bf 127} (2007), 83-102.

\bibitem{Y15} J.-H. Yang, {\em A partial Cayley transform for Siegel-Jacobi disk,}
J. Korean Math. Soc. {\bf 45}, No. 3 (2008), 781-794 $\textit{or}$
arXiv:math.NT/0507216.


\bibitem{Y16} J.-H. Yang, {\em Invariant metrics and Laplacians on Siegel-Jacobi disk,}
arXiv:math.NT/0507217 v1 $\textit{or}$ Chinese Annals of Mathematics, Ser. B,
DOI:10.1007/s11401-008-0348-7, Springer-Verlag Berlin Heidelberg (2009).

\bibitem{Y17} J.-H. Yang, {\em Theta Series Associated With the Schr{\"o}dinger-Weil
Representation,} arXiv:0709.0071v2 [math.NT] 3 Aug 2009.

\bibitem{Y18} J.-H. Yang, {\em The Schr{\"o}dinger-Weil
Representation and Jacobi Forms of Half-Integral Weight,} arXiv:0908.0252v1 [math.NT]
3 Aug 2009.



\bibitem{Yu} H. P. Yuen, Phys. Rev. A, {\bf 13} (1976), 2226--2243.


\bibitem{Zi} C. Ziegler, {\em Jacobi Forms of Higher
Degree}, Abh. Math. Sem. Hamburg $ \textbf{59}$ (1989), 191--224.

\end{thebibliography}

%\vskip 1cm \noindent $ \textsf{Department of Mathematics}$ \\ $ \textsf{Inha University}$ \\  $ \textsf{Incheon 402-751}$\\
%$ \textsf{Republic of Korea}$ \vskip 0.31cm \noindent

%\noindent $   \texttt{email\, :\,jhyang@inha.ac.kr}$

\end{document}